\RequirePackage{fix-cm}
\documentclass[smallcondensed]{svjour3}     
\smartqed  
\usepackage{ulem,amsmath, amsfonts,amssymb,graphicx,color,rotating,url,tikz, pgffor, xcolor,xkeyval,hyperref, listings, mathtools, etoolbox,bm, aliascnt,epstopdf, float,algorithm, booktabs}
\usepackage[noend]{algpseudocode}
\usepackage[numbers]{natbib}

\usepackage[nolist]{acronym}
\begin{acronym}
\acro{AC}{Attainability Condition}
\end{acronym}

\begin{document}

\title{Inventory Control with Modulated Demand and a Partially Observed Modulation Process
}

\titlerunning{Inventory control with Markov-modulation process}        

\author{Satya S. Malladi       \and
        Alan L. Erera          \and
        Chelsea C. White III 
}


\institute{S.S. Malladi \at
              Kantar Analytics Practice, India\\
              Tel.: +91 79896 25293\\
              \email{sarvanisatya@gmail.com}           
           \and 
           A.L. Erera \at
           School of Industrial and Systems Engineering, Georgia Institute of Technology,  USA \\
           \and
           C.C. White III  \at
           School of Industrial and Systems Engineering, Georgia Institute of Technology,  USA
}


\maketitle

\begin{abstract}
We consider a periodic review inventory control problem having an underlying modulation process that affects demand and that is partially observed by the uncensored demand process and a novel additional observation data (AOD) process. We present an attainability condition, AC, that guarantees the existence of an optimal myopic base stock policy if the reorder cost $K=0$ and the existence of an optimal $(s, S)$ policy if $K>0$, where both policies depend on the belief function of the modulation process. Assuming AC holds, we show that (i) when $K=0$, the value of the optimal base stock level is constant within regions of the belief space and that each region can be described by two linear inequalities and (ii) when $K>0$, the values of $s$ and $S$ and upper and lower bounds on these values are constant within regions of the belief space and that these regions can be described by a finite set of linear inequalities. A heuristic and bounds for the $K=0$ case are presented when AC does not hold. Special cases of this inventory control problem include problems considered in the Markov-modulated demand and Bayesian updating literatures.

\end{abstract}
\DeclarePairedDelimiterX{\abs}[1]{\lvert}{\rvert}{\ifblank{#1}{{}\cdot{}}{#1}}
\section{Introduction}
We consider a periodic review, data driven inventory control problem over finite and infinite planning horizons with instantaneous replenishment. We assume that there are several interconnected processes: the completely observed inventory process that keeps track of the inventory level, the uncensored demand process, the action process that represents replenishment decisions, the underlying modulation process that affects demand, and the additional observation data (AOD) process that together with the demand process partially observes the modulation process. The inventory, demand, and action processes are common to
inventory control problems. 

When completely observed by the demand and AOD processes, the modulation process models the case where demand is Markov-modulated. When the modulation process is only observed by the demand process and is assumed static, then the model conforms to a model considered by the Bayesian updating literature. Thus, the model presented in this paper generalizes models found in both the Markov modulated literature and the Bayesian updating literature and hence serves as a bridge between these two major research directions. 

Additionally, the new model can consider scenarios that have not been considered thus far in the literature. For example, assume the modulation process is partially observed by both the demand and AOD processes and is dynamic, which is consistent with reality when the modulation process models the macro economy and the observations of the AOD process are housing starts, consumer spending, and other indicators that partially observe the state of the macro economy.

Noting that improving demand forecasting is of major on-going interest to industry, the modulation and AOD processes can serve as a modeling basis for improved data-driven demand forecasting using data that includes past demand data. The modulation process can represent an unknown static parameter or index of the demand process and dynamic exogenous and partially observed factors, such as the weather, seasonal effects, and the underlying economy. The AOD process can model observations of the modulation process other than demand; e.g., weather, season, day of the week, and macro-economic indicators, and is motivated by the fact that supply chains are becoming increasingly data driven to support improved real-time supply chain control.  For example, consider the demand for housing construction materials.  This demand depends directly on the number of new housing starts.  Demand for new housing is influenced by the state of the economy, and the state of the economy can be inferred from various macro-economic indicators, such as the rate of growth of the Gross Domestic Product, interest rates, and measures of consumer confidence.   Historical economic data can be used to determine the parameter values of the modulation and AOD processes.  

We now outline and present the contributions found in the remainder of this paper.  In \autoref{sec2}, we model this inventory control problem as a partially observed Markov decision process (POMDP).  We show that for the single period problem, there exists an optimal base stock policy, the value of the optimal base stock level is dependent on the belief function but constant within regions of the belief space, and these regions can be described by a finite set of linear inequalities. We then present a generalization of the Veinott attainability condition, \acs{AC}.

In \autoref{sec3}, we assume that reorder cost is zero ($K = 0$) and that \acs{AC} holds.  Our first result extends the single period result in \autoref{sec2} to the countable period case and shows that the single period optimal base stock policy is myopic and hence identical for all decision epochs.  We also show that use of this optimal policy and \acs{AC} guarantees that the current base stock level is always at least as great as the current inventory level.  We define a partial order on the belief space and present conditions that guarantee the base stock level is non-decreasing with regard to the partial order.   

In \autoref{sec4}, we assume $K=0$ and that \acs{AC} does not hold.  Preliminary analysis shows that there is an optimal base stock policy for this case; however, this base stock policy is not myopic and is hence stage dependent. This complication identifies an interesting future research challenge and motivates the search for a good, easily computed, and easily implemented heuristic. The heuristic we chose to investigate is the easy to determine and easy to implement optimal policy determined in \autoref{sec3}. We determine a lower bound on the optimal value function and use the policy determined in \autoref{sec3} to generate an upper bound on the optimal value function. We then present an upper bound on the difference between the upper and lower bounds on the optimal value function. 
We also show that the upper bound on the optimal cost function is piecewise linear in the belief function for the finite horizon case but may not be continuous. Hence, and counter intuitively, improved observation quality of the modulation process may not result in improved systems performance using this heuristic.

In an extensive numerical study, we show that \acs{AC} holds for at least 63\% of the instances considered.  Further, when this condition is violated, we show that the difference between the proposed upper and lower bounds on the optimal value function is less than 0.6\% on average and that the heuristic generates an average percentage optimality gap of 0.51\%.  This study supports the claims that (i) \acs{AC} often holds and (ii) when \acs{AC} does not hold, the optimal base stock policy determined in \autoref{sec3} may be an excellent heuristic.

In \autoref{sec5}, we consider the $K > 0$ case and assume throughout that \acs{AC} holds. We show that there exists an optimal $(s,S)$ policy and determine upper and lower bounds on $s$ and $S$ for the finite and infinite horizon cases, where each bound and the values of $s$ and $S$ are dependent on the belief function of the modulation process. Each of these bounds and the values of $s$ and $S$ are shown to be constant within regions of the belief space described by a finite number of linear inequalities. An outline of an approach for determining an optimal $(s,S)$ policy and the resultant expected cost function for the finite horizon case are presented in the e-companion.

Conclusions are presented in \autoref{sec6}. 

\subsection{Literature Review}
This survey is organized around various assumptions made in the literature regarding the modulation process, with emphasis on how the model considered in this paper generalizes much of the Markov-modulated demand and Bayesian updating inventory control literatures. After a brief overview of the general inventory and POMDP literatures, we review the literature concerned with the $K = 0$ case, followed by the $K > 0$ case. 

Inventory control has been studied extensively over six decades; see \cite{graves93,choiNewsvendor, khouja99,petruzzi99,qin2011,bellman58,katehakis2016, huh09, katehakis2012}, and \cite{arrowbook58} for detailed surveys. We also survey several nonparametric approaches. We model the inventory control problem considered in this paper as a POMDP; see \cite{smallwood73} and \cite{sondik78} for the POMDP foundational results on which our results are based.

\subsubsection{The $K=0$ case}  
Assume the modulation process is completely observed and static.  This special case of the problem we consider in this paper, for which assumption \acs{AC} always holds, was first considered by \citet{arrow51}, and \citet{karlin58a}, and \citet{karlin58b}, various  extensions of which are detailed in surveys by \citet{graves93}, \citet{khouja99}, \citet{petruzzi99}, and \citet{qin2011}.

A special case of the problem considered in this paper is obtained by restricting the modulation process to be completely observed and nonstationary with known demand distributions for each period.  This case was first considered by \citet{karlin59a,karlin59b}, \citet{iglehart62}, and \citet{veinott65b, veinott65a}. A base stock policy dependent on the completely observed state of the modulation process (current demand distribution) was proved to be optimal in \cite{karlin59a} and \cite{karlin59b}. We will later see that this result is an implication of the generalized result presented in this paper.  \citet{iglehart62} developed computational approaches for determining the base stock level. \citet{zipkin89} extended these results to the average cost criterion and to cyclic costs. \citet{veinott65a,veinott65b} proved the existence of an optimal myopic base stock policy when the base stock level at the next decision epoch is guaranteed to exceed the current inventory position after satisfying demand (i.e., the attainability assumption) for independent and correlated nonstationary demands across time periods, respectively. \citet{veinott65a} also provided sufficient conditions for this assumption. We generalize Veinott's attainability assumption for the problem considered in this paper (\acs{AC}) and also provide sufficient conditions for this generalized attainability condition. \citet{morton78} studied an inventory system with the additional option of disposal of inventory at a cost and nonstationary demands in each period. 
\citet{lovejoy92} modeled explicit dependence of a generalized demand process on a completely observed modulation process with exogenous parameters and demand history and derived an upper bound on the optimal cost for scenarios such as Markov modulation and additive and multiplicative demand shocks. 
\citet{song93} modeled the modulation process as a completely observed underlying ``state-of-the-world'' in a continuous time framework similar to \cite{iglehart62}, with a Markov-modulated Poisson demand process. A ``state-of-the-world'' dependent base stock policy was proved to be optimal. An optimal myopic policy was shown to exist when the attainability of the next period's base stock level is guaranteed.   \citet{sethi97} extended the results of \cite{song93} to a discrete time system and obtained analogous results.

\citet{dvoretzky52}, \citet{scarf59}, and \citet{murray66} analyzed the case where the modulation process is static, partially observed by the demand process, completely unobserved by the AOD process, and represents unknown parameters of a single stationary distribution. While \citet{scarf59} proved the optimality of a statistic-dependent base stock policy, \citet{murray66} extended the results to determine a Bayesian update on unknown parameters. \citet{azoury84} and \citet{azoury85} extended these results (\cite{scarf59}) to other distributions and compared this method with non-Bayesian mixture methods.  \citet{lovejoy90} proved the optimality of a myopic base stock policy for ``parameter adaptive models'' of demand,  and \citet{lariviere99} dealt with an unknown stationary distribution of demand partially observed by a scale parameter and a shape parameter. \citet{kamath2002} presented a study of the Bayesian updating mechanism  with and without nonstationarity. We note that when the modulation process is static, although partially observed, the attainability assumption always holds. We remark that the problem instances of the problem we consider in this paper, for which \acs{AC} is sometimes violated, always involve a dynamic modulation process. 

Partial observability of demand outcomes, which is different from partial observability of modulation process, results from limitations on the accuracy of inventory book-keeping (in \cite{benso07b} and \cite{ortiz13}), and censoring (in \cite{benso07a, benso08,huh09a, huh11, katehakis2015, zhang2020, yuan2021}).  \citet{benso07a, benso08} treated Markovian modulation of demand as a special case. Their problem formulation differs from our framework in that their demand process (not the modulation process) is partially observed (censored). \citet{ding2002} presented an analysis of optimal policies for the Bayesian newsvendor problem with and without censoring. 

For the case where the modulation process is partially observed by the demand process, completely unobserved by the AOD process, and dynamic, \citet{treharne2002} proved the existence of an optimal state-dependent base stock policy for an un-capacitated inventory system.
\citet{arifoglu2010}  proved the optimality of inflated state-dependent base stock policies for capacitated production systems under Markov-modulated demand and supply processes (extending \cite{gallego2004}). \citet{bayraktar2010} studied a completely unobserved Markov-modulated Poisson demand process in a continuous-review inventory system with reorder cost and lost sales (censoring). 
 
\citet{treharne2002} and \citet{arifoglu2010}, however, did not prove the existence of an optimal \textit{myopic} state-dependent base stock policy, which we prove in this paper, assuming \acs{AC} holds.  Further, we show that the belief space can be partitioned into subsets by a finite set of linear inequalities and that the base stock level is constant within each of these subsets. Such regions have also been observed in the numerical example provided in \cite{bayraktar2010}; however, no explanation is given for such behavior. The linear partition of the belief space we present provides an easily computed approach to determine an optimal base stock level for any given belief vector.

\subsubsection{The $K>0$ case} 
\citet{scarf60} and \citet{iglehart63} proved that there exists an optimal $(s, S)$ policy under finite and infinite horizons, respectively. \citet{iglehart63} presented the first set of bounds on period-wise reorder points and base stock levels, which were later tightened by \citet{veinott65c}. \citet{veinott65d} extended \citet{veinott65a, veinott65b} to the $K>0$ case. More recently, \citet{chenroundy15} presented sufficiency conditions of divergence and $K$-convexity for the optimality of $(s, S)$ policies under time-varying parameters and correlated demand variables modulated by an underlying ``state-of-the-world'' variable. Our results extended to the $K>0$ case lead to significantly reduced computational effort in determining the optimal policy compared to \citet{bayraktar2010} when \acs{AC} holds. 

More recently, data-driven and nonparametric approaches for describing demand uncertainty have garnered interest. Related literature includes research presented in the following papers:  \cite{bookbinder89,gallego93,godfrey01,perakis2008,huh11,besbes13,klabjan13,ban15}, \cite{bertsimas14,ban14,levi15,cheung15,xin2015,mamani16,ferreira16}. 
Future research may involve a blend of nonparametric approaches with Bayesian approaches, such as the data driven approach presented in this paper. 

\section{Problem Description and Preliminary Results \label{sec2}}
We describe the inventory control problem in \autoref{probdesc}. We then model the problem as a POMDP and present optimality equations and other standard results in \autoref{model}. In \autoref{lxy} we present results associated with the single period expected cost function that will be useful in later sections and also present the condition \acs{AC}. 
\subsection{Problem Definition\label{probdesc}}
We consider an inventory control problem that involves the inventory process $\{s(t), t=0,1,\dots\}$, the modulation process $\{\mu(t), t=0,1,\dots\}$, the demand process $\{d(t), t=1,2,\dots\}$, the additional observation data (AOD) process  $\{z(t), t=1,2,\dots\}$, and the action process $\{a(t), t= 0,1,\dots\}$. These processes are linked by the state dynamics equation
 $s(t+1) = f\big( y(t), d(t+1)\big) $, where $y(t) = s(t)+a(t)$, and the given conditional probability $\text{Pr}\big(d(t+1),z(t+1),\mu(t+1) \mid \mu(t) \big)$. 
We assume the single period cost accrued between decision epoch $t$ and $t+1$ is $c\big(y(t), d(t+1)  \big)$, where 
$c(y, d)$ is convex in $y$ and $\lim_{\abs{y}\to \infty} c(y,d) \to \infty$ for all $d$. We also assume that $c(y,d)$ is piecewise linear in $y$ for all $d$ and that the facets describing $c(y,d)$ intersect at integers. We will have particular interest in the case where $f(y,d) = y-d$, which assumes backlogging, and $c(y,d) = p\big(d-y\big)^+ + h\big(y-d\big)^+,$ where 
 $p$ is the shortage penalty per period for each unit of stockout, 
$h$ is the holding cost per period for each unit of excess inventory after demand realization, and $(g)^+ = \max(g,0)$. 
Without loss of significant generality, this definition of single period cost does not include an ordering cost.  It is straightforward to transform an inventory problem with a strictly positive ordering cost into an inventory problem with no ordering cost for a wide variety of cost and dynamic models of inventory position, e.g., $ f(y, d) = y - d$ or $f(y, d) = (y- d)^+$ and $c(s, y, d) = c'(y - s) + p(d -y)^+ + h(y - d)^+$, where in this case the single period cost accrued between decision epochs is dependent on $s$ and $c'$ is the cost per unit ordered.  

We assume that the modulation, demand and AOD state spaces are all finite, the inventory process has a countable state space, and the action space is the set of non-negative integers. We assume the action at $t$ can be selected based on $s(t),d(t),d(t-1),\dots,z(t),z(t-1),\dots$, and the prior probability mass vector $\{\text{Pr}\big(\mu(0)=\mu_i\big),\forall i\}$. Thus, the inventory process is completely observed, demand is not censored, and the modulation process is partially observed by the demand and AOD processes. The problem is to determine a policy that minimizes the expected total discounted cost over the infinite horizon, where we let $\beta \in [0,1)$ be the discount factor. It is assumed throughout that replenishment is instantaneous. 

We remark that the inventory, demand, and action processes are all part of inventory control problems considered in the literature. As indicated in the literature review, the modulation process is also part of the structure of inventory control problems with Markov-modulated demand. The AOD process is intended to provide information about the modulation process, where appropriate, in addition to that provided by the demand process, such as macro-economic data. Throughout we assume demand realization is uncensored and completely revealed.  This assumption is in contrast to the censored demand case where only sales data are available to the decision maker. 

We note that the conditional probability $\text{Pr}\big(d(t+1),z(t+1),\mu(t+1) \mid \mu(t) \big)$ is the product of two conditional probabilities:
\begin{enumerate}
\item $\text{Pr}\big(d(t+1),z(t+1) \mid \mu(t+1),\mu(t)\big)$, the demand and AOD probabilities, conditioned on the modulation process
\item $\text{Pr}\big(\mu(t+1)\mid \mu(t)\big)$, the state transition probabilities for the (Markov-modulated) modulation process.
\end{enumerate}
The Baum-Welch algorithm is typically used to estimate parameters of a POMDP, \textit{viz.}, observation and transition probabilities and initial belief state (see \cite{atrash2010} for a review on POMDP training methods). 
 
\subsection{The POMDP Model and Preliminary Results\label{model}}
This problem can be recast as a partially observed Markov decision problem as follows. Results in \cite{smallwood73} and \cite{sondik78} imply that $(s(t),x(t))$ is a sufficient statistic, where $N$ is the number of values the modulation process can take, the belief function $\bm{x}(t) = \text{row} \{x_1(t), \dots, x_N(t) \}$, is such that $x_i(t) = \text{Pr}\big(\mu(t) = \mu_i \mid d(t), \dots, d(1), z(t), \dots, z(1), x(0) \big)$, and $\bm{x}(t) \in X = \{\bm{x}\in \mathbb{R}^N: \bm{x}\geq 0 \text{ and } \sum_{i=1}^N x_i = 1 \}$. For $\bm{g}\in \mathbb{R}^N$, let $\bm{g}\underline{1} = \sum_{n=1}^N g_n$. Thus, the inventory process is completely observed, the modulation process is partially observed through the demand and AOD processes, and the state of the modulation process is characterized by the belief function.  Let
\begin{eqnarray}
\nonumber {P}_{ij}(d,z) &=& \text{Pr}\big(d(t+1)=d, z(t+1)=z, \mu(t+1) = j \mid \mu(t) = i \big), 
\\ \nonumber  \bm{P}(d,z) &=& \{{P}_{ij}(d,z) \},
\\ \nonumber \sigma(d,z,\bm{x}) &=& \bm{x}\bm{P}(d,z) \underline{1} = \sum_{i=1}^N x_i \sum_{j} {P}_{ij}(d,z), 
\\ \nonumber \bm{\lambda}(d,z,\bm{x}) &=& \text{row} \big\{ \lambda_1(d,z,\bm{x}), \dots, \lambda_N(d,z,\bm{x}) \big \} = \bm{x}\bm{P}(d,z)/\sigma(d,z,\bm{x}), \sigma(d,z,\bm{x}) \neq 0,
\\ \nonumber L(\bm{x},y) &=& E \big[c(y,d) \big]=\sum_{d,z}\sigma(d,z,\bm{x}) c(y,d).
\end{eqnarray}
\vspace{-2ex}
Define the operator $H$ as 
\begin{eqnarray}
\label{eq1} [Hv](\bm{x},s) = \min_{y\geq s} \bigg \{ L(\bm{x},y) + \beta \sum_{d,z} \sigma (d,z,\bm{x}) v(\bm{\lambda}(d,z,\bm{x}),f(y,d)\big) \bigg \}. 
\end{eqnarray}
Results in \cite{puterman94} guarantee that there exists a unique cost function $v^*$ such that $v^* = Hv^*$ and that this fixed point is the expected total discounted cost accrued by an optimal policy. We can restrict search for an optimal policy to t-invariant functions that select $a(t)$ on the basis of $\big(s(t),\bm{x}(t)\big)$, the function $\psi$ such that $\psi \big(s(t),\bm{x}(t) \big) = a(t)$ causing the minimum in \eqref{eq1} to be attained is an optimal policy, and $\lim_{n\to \infty} \| v^* - v_n \| = 0$, where the (finite horizon) cost function $v_{n+1} = Hv_n$ for any given bounded function $v_0$ and $\|.\|$ is the sup-norm. The function $L(\bm{x}, y)$ is the expected single period cost, conditioned on belief $\bm{x}$ and inventory level $y$.  From the perspective of Bayes' Rule, note that $\bm{x} = \bm{x}(t)$ can be thought of as the prior probability mass function of $\mu(t)$, $\sigma(d, z, \bm{x})$ is the probability that the demand and AOD processes will have realizations $d = d(t+1)$ and $z = z(t+1)$, respectively, given $\bm{x}$, and $\bm{x}(t+1) = \bm{\lambda}(d, z, \bm{x})$ is the posterior probability mass function of $\mu(t)$, given $d$, $z$, and $\bm{x}$.

With respect to optimal value function structure, results in  \cite{smallwood73} guarantee that $v_n(\bm{x},s)$ is piecewise linear and concave in $\bm{x}$ for each fixed $s$ for all finite $n$, assuming $v_0(\bm{x},s)$ is also piecewise linear and concave in $\bm{x}$ for each $s$. In the limit $v^*(\bm{x},s)$ may no longer be piecewise linear in $\bm{x}$ for each $s$; however, concavity will be preserved. 
 
Regarding the value of information, results in \citep{sondik78} and \citep{chang15b} can be used to determine upper and lower bounds on $v_n(\bm{x},s)$ for all $n$ and $v^*(\bm{x},s)$ as a function of observation quality and that as observation quality improves, optimal expected systems performance will never degrade. This result may not hold for some sub-optimal policies, as discussed in \autoref{eg4} of \autoref{sec_vu}. 
\subsection{$L(x,y)$ Analysis for Backordering Systems \label{lxy}}
We now examine $L(\bm{x},y)$ in more detail for backordering systems where $f(y,d) = y-d$ and $c(y,d) = p(d-y)^++h(y-d)^+$. Let $\{d_1, \dots, d_M\}$ be the set of all possible demand values, where $d_m < d_{m+1}$, for all $m = 1,\dots, M-1$.  Letting $\sigma(d,\bm{x}) = \sum_z \sigma(d,z,\bm{x})$, define for all $m=0, \dots, M$,
\begin{eqnarray*}
A_m(\bm{x}) &=& h\sum_{k=1}^m \sigma(d_k,\bm{x}) - p\sum_{k=m+1}^M\sigma(d_k,\bm{x}),
\\B_m(\bm{x}) &=& p\sum_{k=m+1}^Md_k\sigma(d_k,\bm{x}) - h \sum_{k=1}^m d_k \sigma(d_k,\bm{x}).
\end{eqnarray*}
Note, $A_0(x) = -p$ and $B_0(x) = p\sum_{k=1}^M d_k\sigma(d_k,\bm{x})$. Proof of the next result, which provides structure that will prove useful, is straightforward. 
\begin{lemma}\label{lemma1}
For all $\bm{x}\in X$: 
\begin{description}
\item[(i)] $
L(\bm{x},y) = \begin{cases} A_0(\bm{x})y+B_0(\bm{x}) =  p\sum_{k=1}^M\sigma(d_k,\bm{x}) (d_k - y), & y\leq d_1
\\A_m(\bm{x})y+B_m(\bm{x}), & d_m\leq y \leq d_{m+1}, 
\\ & m = 1, \dots, M-1
\\A_M(\bm{x})y+B_M(\bm{x}) =  h\sum_{k=1}^M\sigma(d_k,\bm{x}) (y-d_k), & d_M \leq y.
\end{cases}
$
\item[(ii)] for all $m = 1, \dots, M-1$, $A_{m+1}(\bm{x}) = A_m(\bm{x}) + (h+p)\sigma(d_{m+1},\bm{x})$, 
and hence, $A_{m+1}(\bm{x}) \geq A_m(\bm{x})$.
\item[(iii)] for all $m = 1,\dots, M-1$, $B_{m+1}(\bm{x}) = B_m(\bm{x}) - (p+h) d_{m+1} \sigma(d_{m+1},\bm{x}),$ 
and hence, $B_{m+1}(\bm{x}) \leq B_m(\bm{x})$.
\item[(iv)] for all $m = 1, \dots, M$, $A_{m-1}(\bm{x})d_m + B_{m-1}(\bm{x}) = A_m(\bm{x})d_m + B_m(\bm{x})$. 
\item[(v)] $L(\bm{x},y) = \max_{0\leq m \leq M} \bigg[ A_m(\bm{x})y+B_m(\bm{x})\bigg].$
\end{description}
\end{lemma}
\subsubsection{\label{lin_part} Base Stock Policy: Linear Partition of Belief Space.} 
\autoref{lemma1} establishes that $L(\bm{x},y)$ is piecewise linear and convex in $y$ for all $\bm{x}\in X$. Let $s^*(\bm{x})$ be the smallest integer that minimizes $L(\bm{x},y)$ with respect to $y$. Note that it is sufficient to restrict $s^*(\bm{x})$ to the set $\{d_1, \dots, d_M\}$, i.e., $\{s^*(\bm{x}): \bm{x} \in X \} \subseteq \{d_1, \dots, d_M \} $. Hence, 
$L(\bm{x},s^*(\bm{x})) = \min_{1\leq m \leq M} \bigg\{A_m(\bm{x})d_m + B_m(\bm{x}) \bigg \}$. 
Let $\mathcal{P}_1$ be the partition of $X$ composed of elements $X_m = \{\bm{x}\in X: s^*(\bm{x})=d_m\}$. Thus, $\mathcal{P}_1 = \{X_m, m=1,\dots,M\}$, where $X_m$ is non-null for all $d_m$ such that $\min\{s^*(\bm{x}): \bm{x}\in X\}\leq d_m\leq \max\{s^*(\bm{x}): \bm{x}\in X\}$. We characterize $X_m$ as follows. 
\begin{lemma}\label{lemma2} Let $\bm{P}(d) = \sum_z \bm{P}(d,z), \ \forall \ d$. 
For $m=1,\dots, M$, 
\begin{eqnarray}
\label{newsv} X_m = \bigg \{ 
\bm{x} \in X: \bm{x}\sum_{k=1}^{m-1} \bm{P}(d_k)\underline{1} < p/(p+h) \leq \bm{x}\sum_{k=1}^m \bm{P}(d_k) \underline{1} \bigg \}.
\end{eqnarray}
\end{lemma}
Note that the criterion in \eqref{newsv} can be re-written as:
$$\sum_{k=1}^{m-1}\sigma(d_k, \bm{x}) < p/(p+h) \leq \sum_{k=1}^m\sigma(d_k,\bm{x}),$$ 
where $\sigma(d_k, \bm{x})$ is the probability of observing demand outcome $d_k$ when the current belief is $\bm{x}$. This criterion is identical to the \textit{newsvendor problem}'s criterion for determining the optimal base stock policy with the probability mass function of demand given by $\sigma(d_k,x), \ \forall \ k$. 

Due to the linearity of $\sigma(d, \bm{x})$ in $\bm{x}$, the above criterion results in a \textit{linear partition} of the belief space. We note that the partition thus obtained is independent of the values of demand and AOD outcomes but depends only on the parameters, $P_{ij}(d,z)$, $p$, and $h$. We remark that $X_m$ for all $m$ can be described by two inequalities linear in $\bm{x}$, which is true irrespective of the values $N$ and $M$ take, since $L(\bm{x},y)$ is piecewise linear in $\bm{x}$ for fixed $y$.  

\begin{example} \label{eg1} Let $ M=7, N=3, h= 1, p = 3$, $\bm{d}=[5, 10, 15, 20, 25, 30, 35]$ and 
\begin{align*}
& \bm{P}=\begin{bmatrix}
0.0192	&0.8744 &	0.1063 \\
0.0437	&0.4712 &	0.4851 \\ 
0.4467	&0.0313&	0.522 
\end{bmatrix},
\\ & \bm{Q^D} = \begin{bmatrix}
0.0207&	0.2321&	0.0717&	0.2054&	0.1519&	0.0346&	0.2837
\\0.2697	&0.208	&0.2044&	0.1942&	0.0748&	0.0427&	0.0062
\\0.0283&	0.0378&	0.0429&	0.0605&	0.1335&	0.3001&	0.3969
\end{bmatrix}
\\ & \text{where } \bm{P} = \{P_{ij}\} \text{ and } P_{ij} = \text{Pr}\big(\mu(t+1) = j \mid \mu(t) = i\big)
, 
\\ & \bm{Q^D} = \{q_{jd}^D \}, \  q_{jd}^D = \text{Pr}\big(d(t+1) = d \mid \mu(t+1) = j\big), 
\text{ where } q_{jd}\text{ is independent of } i.  
\end{align*}
Note that $P_{ij}(d) = P_{ij} q^D_{jd}$ is independent of $z$. The belief space is given by the triangle with vertices $(1,0,0)$, $(0,1,0)$, and $(0,0,1)$ (described by $x_1+x_2+x_3=1$, $x_1 \geq 0$, $x_2 \geq 0$, and $x_3\geq 0$), where modulation state $n+1$ indicates a stronger economy than modulation state $n$, for all $n$.  
\begin{figure}[h]
\centering
\includegraphics[scale = 0.6, trim={65pt 65pt 0pt 0pt},clip]{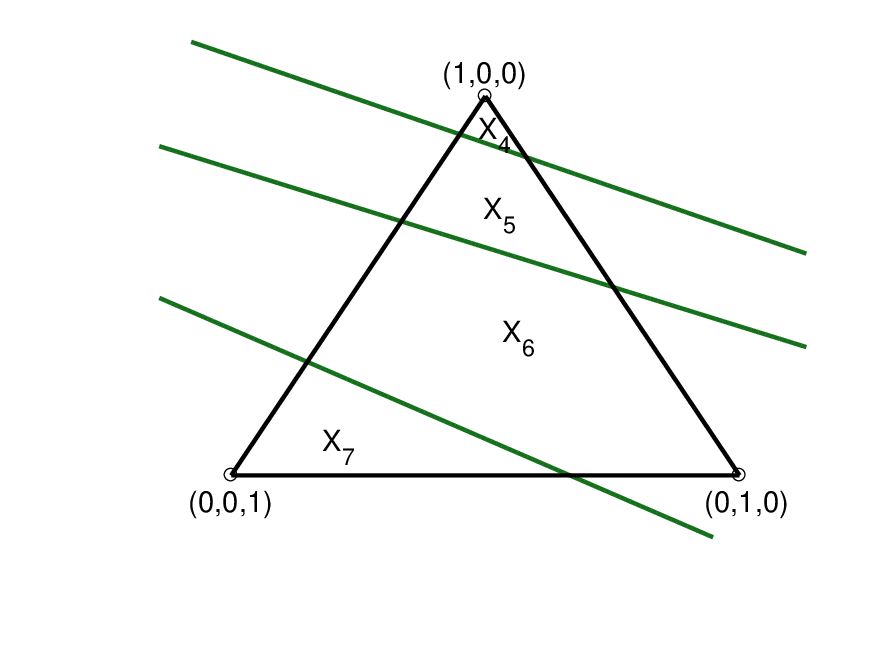}
\caption{Example of $\mathcal{P}_1$ with $N=3$ and $M=7$}
\label{fig1} 
\end{figure}
\autoref{fig1} depicts the belief space, $X$, overlaid with the  partition, $\mathcal{P}_1$ (derived in \autoref{lemma2}).  $\mathcal{P}_1$ divides $X$ into 4 regions of constant base stock level, \textit{viz}., $X_4$ through $X_7$. For any belief vector in $X_m$, the optimal order-up-to level for the one period problem is $d_m$. Hence, the optimal myopic base stock levels are $20$, $25$, $30$, and $35$ in $X_4$, $X_5$, $X_6$, and $X_7$ respectively. 

If the AOD process is dependent on $\mu(t+1)$ (e.g. current state of the economy), and has two outcomes $\{z_1, z_2\}$ (e.g. real estate price levels) with $\bm{R^Z} = \{r_{jz}^Z\}$, $r_{jz}^Z =\text{Pr} \big(z(t+1)=z \mid \mu(t+1) = j\big)$, then $P_{ij}(d,z) = P_{ij} \ q_{jd}^D \ r_{jz}^Z$. 
The regions presented in \autoref{fig1} do not change as \eqref{newsv} does not depend on the values of the outcomes of either the AOD or demand processes. Note $\lambda(d_1, e_3) = [0.28, 0.26, 0.46]$ when the AOD process is uninformative. Let $\bm{R^Z}= [1,  0;\ 0.7, 0.3; \ 0, 1]$ and hence the AOD process is informative. Then, $ \lambda(d_1, z_1, e_3) = [0.61, 0.39, 0 ]$ and $\lambda(d_1, z_2, e_3) = [0, 0.15, 0.85] $. We note that the availability of additional observation data leads to substantially different updated belief functions.
\end{example}

We now present an assumption that would help establish the optimality of the myopic base stock policy in \autoref{sec3}.
\subsection{Definition of the Attainability Condition (AC)}
We now present the attainability condition, \acs{AC}. Let $s^*(\bm{x})$ define an order-up-to level $y$ given a belief state $\bm{x}$.
\begin{att_con}[\textbf{AC}]
\label{A1} $f\big(s^*(\bm{x}),d\big) \leq s^*\big(\bm{\lambda}(d,z,\bm{x})\big)$ for all $\bm{x}\in X$, $d \in \{d_1, \dots, d_M\}$, and $z$. 
\end{att_con}
\acs{AC} assumes that the amount of inventory after demand is satisfied, which is $f\big(s^*(\bm{x}),d\big)$, never exceeds the order-up-to level  at the next decision epoch, which is $s^*(\bm{\lambda}(d,z,\bm{x}))$. This assumption is  identical to assumptions made in the inventory literature specialized to the problems considered in \cite{veinott65b, veinott65a}. 
\section{When \acs{AC} holds}\label{sec3}
\subsection{Main Result}
We assume throughout this section that \acs{AC} holds. We now present the main result of this section. Proof of the following result is provided in the e-companion.  
\begin{proposition}\label{prop1}
Assume \acs{AC} holds, $L(\bm{x},y)$ is piecewise linear in $y$ for all $\bm{x}\in X$, $s^*(\bm{x})$ is the smallest integer that minimizes $L(\bm{x},y)$ with respect to $y$, and that $f(y,d)$ is non-decreasing in $y$ for each $d$. Then, $v_n(\bm{x},s) = v_n\big(\bm{x},\max\{s^*(\bm{x}),s\}\big)$ is non-decreasing and convex in $s$ for all $n$ and $\bm{x}$. Further, the myopic base stock policy that orders up to $\max\{s^*(\bm{x}),s\}$ is an optimal policy. 
\end{proposition}
Thus, when \acs{AC} is satisfied and recalling that $s^*(\bm{x})$ is determined using the inequalities presented in (\ref{newsv}), ordering up to $\max\{s^*(\bm{x}), s\}$  at every decision epoch is optimal for any finite horizon problem and the infinite horizon problem. This result ensures  significantly less computational effort for computing the optimal policy compared to the procedure proposed in \cite{bayraktar2010} for the special case where the modulation process is completely unobserved by the AOD process. 
\subsection{\acs{AC} Analysis}
Proof of the following preliminary result is straightforward. 
\begin{lemma}
Assume \acs{AC} holds, apply the base stock policy ``order up to $\max\{s^*(\bm{x}),s\}$'', and assume $s(t) \leq s^*(\bm{x}(t))$. Then, $s(\tau) \leq s^*(\bm{x}(\tau))$ for all $\tau \geq t$.  \label{lemma3}
\end{lemma}
Thus, once the inventory level falls at or below the base stock level, \acs{AC} guarantees that the inventory level will always fall at or below the base stock level at the next decision epoch.

In \autoref{eg2}, we show that \autoref{eg1} satisfies \acs{AC}, following a  preliminary result. Define the binary operator for first order stochastic dominance, $\preceq$, as follows: for $\bm{x}, \bm{x'}\in X$, $ \bm{x} \preceq \bm{x'} \iff \sum_{i=n}^N x_i \leq \sum_{i=n}^N x_i'  \ \ \ \forall \ n=1,\dots,N.$

\begin{lemma}
\label{lemma6}  Assume if $i\leq i'$, then $ \sum_{k=m}^M \sum_j P_{ij}(d_k) \leq \sum_{k=m}^M \sum_j P_{i'j}(d_k) $ for all $m=1,\dots,M$. Then, $\bm{x}\preceq \bm{x'}$, implies $s^*(\bm{x}) \leq s^*(\bm{x'})$.  
\end{lemma}

\autoref{lemma6} guarantees that $s^*(\bm{e}_N) \geq s^*(\bm{x})$ for all  $ \bm{x} \in X \}$, where $\bm{e}_n \in X $ has a 1 as its $n^{\text{th}}$ entry. Let $\widehat{\bm{x}}^{d,z} \in X$ be such that $\widehat{\bm{x}}^{d,z} \preceq \bm{\lambda}(d,z,\bm{x})$ $\forall \  \bm{x}\in X$. A simple linear programming procedure for determining $\widehat{\bm{x}}^{d,z}$ is presented in \autoref{app_sec3} of the Appendix.
We can now show that the problem presented in \autoref{eg1} satisfies \acs{AC}.

\begin{example}\label{eg2} Consider the problem in \autoref{eg1}. 
It is straightforward to show that the assumption in \autoref{lemma6} is satisfied. If $f(s^*(\bm{e}_N), d) \leq s^*(\widehat{\bm{x}}^{d,z})$ for all $(d,z)$, then \acs{AC} is satisfied since  
$$ f(s^*(\bm{x}), d) \leq f(s^*(\bm{e}_N), d) \leq s^*(\widehat{\bm{x}}^{d,z})  \leq s^*(\bm{\lambda}(d,z,\bm{x})).$$ 
Use of the procedure in \autoref{app_sec3} of the Appendix verifies that $f(s^*(\bm{e}_N), d) \leq s^*(\widehat{\bm{x}}^{d,z}) \ \forall \ (d, z)  $, and hence, \acs{AC} is satisfied. 
 \autoref{fig2} shows $\{ \bm{\lambda}(d_5,z,\bm{x}): \bm{x}\in X\}$ and $\widehat{\bm{x}}^{d_5, z}$, where both $\widehat{\bm{x}}^{d,z}$ and $\bm{\lambda}(d,z,\bm{x})$ are independent of $z$ (e.g., where $\text{Pr}\big( z(t+1) \mid \mu(t+1), \mu(t) \big)$ is independent of $\mu(t+1)$ and $\mu(t)$). 
 \begin{figure}[h]
\centering
\includegraphics[scale = 0.6, trim={65pt 65pt 0pt 0pt},clip]{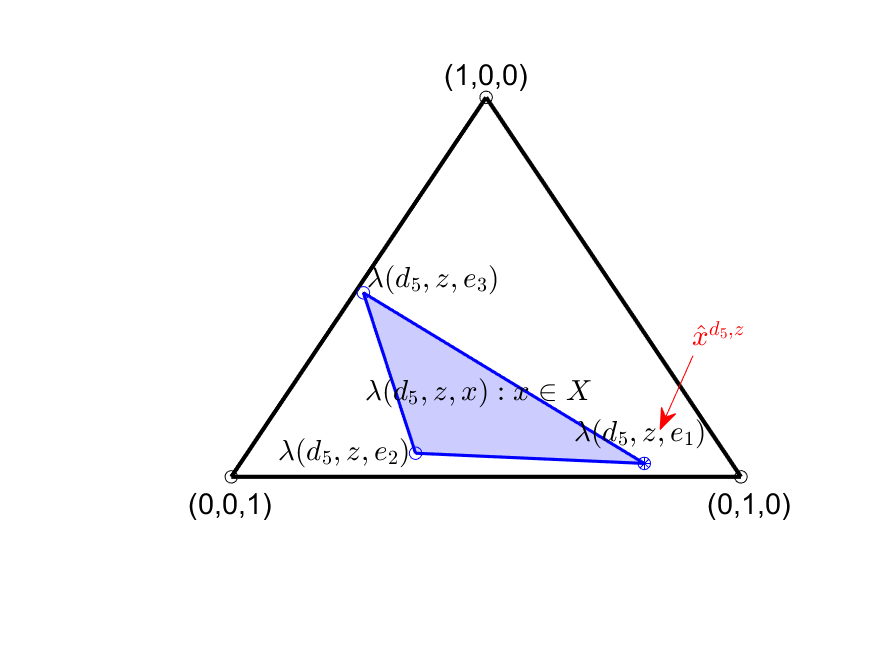}
\caption{$\widehat{\bm{x}}^{d_5, z}$ and $\{\bm{\lambda}(d_5,z,\bm{x}): \bm{x}\in X\}$}
\label{fig2}
\end{figure}
 Thus, since Example 1 satisfies \acs{AC}, ordering up to the myopic base stock level
 at every decision epoch is optimal over finite and infinite horizons.
\end{example}

We remark that for the case where the modulation process is completely observed and $f(y,d)$ is non-increasing in $d$ for all $y$, \acs{AC} is equivalent to $f\big(s^*\big(\mu(t)\big),d(t+1) \big) \leq s^*\big(\mu(t+1)\big)$ for all $d(t+1)$, and hence $f\big(s^*\big(\mu(t)\big), d_1\big) \leq s^*\big(\mu(t+1)\big)$. Note that this is equivalent to the attainability assumption presented by \citet{veinott65b, veinott65a} that guarantees the optimality of a myopic base stock policy for the completely observed nonstationary case. 

A computational procedure to determine the value function of the POMDP is presented in \autoref{app_sec3} of the Appendix. 
\section{When \acs{AC} Does Not Hold\label{sec4}}
We now consider the case where \acs{AC} does not hold. Preliminary analysis shows that there is an optimal base stock policy for this case; however, this base stock policy is not myopic and is hence stage dependent. This complication identifies an interesting future research challenge and motivates the search for a good, easily computed and easily implemented heuristic. 

We remark that although $f(s^*(\bm{x}), d)\leq s^*(\bm{\lambda}(d,z,\bm{x}))$ may not hold for all $(d,z,\bm{x})$, this inequality holds for all $(d,z)$ for at least a subset of $X$, which supports the conjecture that an optimal policy when \acs{AC} holds may be a good sub-optimal policy when \acs{AC} does not hold. For example, let $X_m \subseteq X$ be such that for $\bm{x}' \in X_m$, $s^*(\bm{x}') \leq s^*(\bm{x}) \ \forall \ \bm{x} \in X$. Then, $f(s^*(\bm{x}), d)\leq s^*(\bm{\lambda}(d,z,\bm{x}))$ for all $(d,z) $ for all $\bm{x} \in X_m$. 

In this section, we address the question: how good is the policy ``order up to $\max\{s^*(\bm{x}), s\}$", which is optimal when \acs{AC} holds, as a heuristic when \acs{AC} does not hold? We begin addressing this question by determining a lower bound $v^L$ on the optimal value function. We use the ``order up to $\max\{s^*(\bm{x}), s\}$" policy to generate an upper bound $v^U$ on the optimal value function. We then present an upper bound on the difference $v^U-v^L$.  Numerical results support the claim that this heuristic is near-optimal for the broad class of problems considered in the numerical analysis.  We also show that the upper bound on the optimal cost function is piecewise linear in the belief function for the finite horizon case but may not be continuous. Hence, and counter intuitively, improved observation quality of the modulation process may not result in improved systems performance using this heuristic.
\subsection{\label{lb1}A Lower Bound, $v^L$}
We now present a lower bound on $v_n(\bm{x},s)$. Let 
$$
[H^Lv](\bm{x},s) = L(\bm{x},s^*(\bm{x})) + \beta \sum_{d,z} \sigma(d,z,\bm{x}) v\big(\bm{\lambda}(d,z,\bm{x}),s^*(\bm{\lambda}(d,z,\bm{x}))\big), $$
 $ \ v_{n+1}^L = H^Lv_n^L, \ v_0^L = 0, $
and $v^L$ be the fixed point of $H^L$, which we note is independent of $s$. 
\begin{proposition}\label{prop2}
For all $\bm{x}, s$, and $n$, $v_n^L(\bm{x}) = v^L_n\big(\bm{x},s^*(\bm{x})\big) \leq v_n(\bm{x},s). $
\end{proposition}
The proof follows from the fact that the controller always brings the inventory to $s^*(\bm{x})$, which is not feasible when the inventory is higher than $s^*(\bm{x})$. We remark that $v_n^L(\bm{x})$ can be computed as was $v_n(\bm{x},s)$ for $s\leq s^*(\bm{x})$, in \autoref{sec3}. A tighter lower bound, dependent on $s$ for $s>s^*(\bm{x})$, would replace $L\big(\bm{x},s^*(\bm{x})\big)$ with $L\big(\bm{x},\max\{s^*(\bm{x}),s\}\big)$ in the definition of the operator $H^L$. However, since such a definition of $H^L$ would complicate later analysis, we have chosen not to use this tighter lower bound in the development of results in the sections to follow.  
\subsection{An Upper Bound, $v^U$ \label{sec_vu}}
\begin{figure}[h]
\centering
\includegraphics[scale = 0.4]{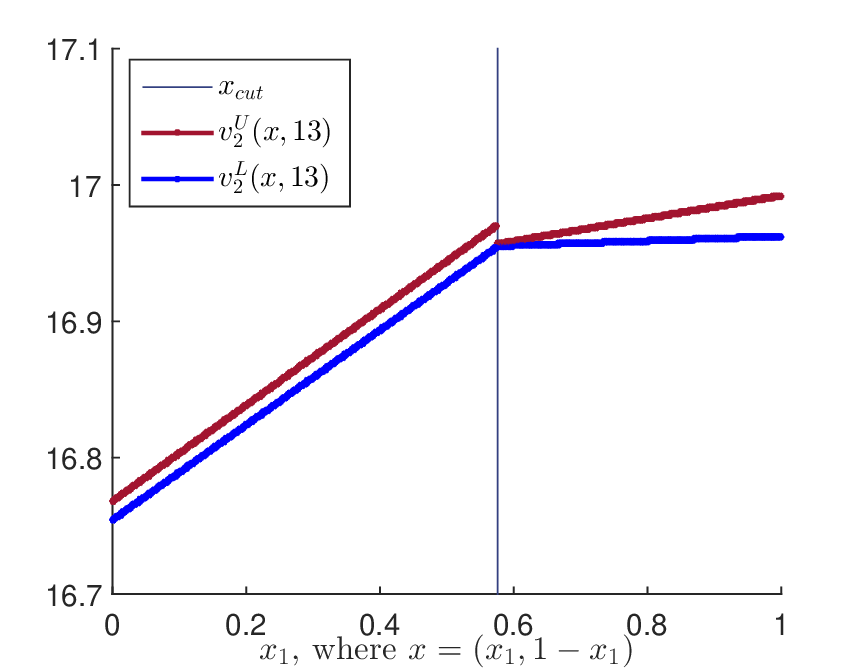}
\caption{ $v^U_2(\bm{x},s)$ and $v^L_2(\bm{x})$ at $s=13$.}
\label{fig4}
\end{figure}

Let
\begin{eqnarray*} [H^Uv](\bm{x},s)&=& L\big(\bm{x},\max\{s^*(\bm{x}),s\} \big) 
\\ &&+ \beta\sum_{d,z} \sigma(d,z,\bm{x}) v\big(\bm{\lambda}(d,z,\bm{x}),f\big(\max\{s^*(\bm{x}),s\}, d \big)\big),
\end{eqnarray*}
 $v^U_{n+1} = H^Uv_n^U, \ v_0^U = 0,$ and let $v^U$ be the fixed point of $H^U$. We remark that $v^U$ is the expected cost to be accrued by the ``order-up-to  $\max\{s^*(\bm{x}),s\}$'' policy, which is feasible but may not be optimal when \acs{AC} is not satisfied, and hence represents an upper bound on the optimal cost function. It is straightforward to prove the following structural result. 
\begin{proposition}
For all $n$ and $\bm{x}$, $v^U_n(\bm{x},s) = v_n^U(\bm{x},s^*(\bm{x}))$ for $s\leq s^*(\bm{x})$, and $v^U_n(\bm{x},s)$ is non-decreasing and convex in $s$. 
\end{proposition}
We prove the following result in the e-companion of this paper. 
\begin{lemma}\label{lemma8.5} For each $n \geq 1$, there is a partition $\mathcal{P}_n$ of $X$ that is defined by a finite set of linear inequalities such that on each element of this partition $v^U_n$ is linear in $\bm{x}$. Further, $\mathcal{P}_{n+1}$ is at least as fine as $\mathcal{P}_n$  (i.e., if $S \in \mathcal{P}_{n+1}$, then there is an $S'\in \mathcal{P}_n$ such that $S\subseteq S'$). 
\end{lemma}
Thus, $v_n(\bm{x},s)$ is piecewise linear in $\bm{x}$ for each $s$. Note that $\mathcal{P}_1$ is defined in \autoref{lin_part}. However, \autoref{eg4} shows that $v_n^U(\bm{x},s)$ may be discontinuous and hence not concave in $\bm{x}$ for each $s$. Thus, according to \citet{white1980}, it may not be true that the improved observation accuracy will improve the performance of the ``order up to $\max\{s^*(\bm{x}),s\}$'' policy if \acs{AC} is not satisfied. 
\begin{example}\label{eg4} Assume $f(y,d) = y-d$, $c(y,d) = p(d-y)^++h(y-d)^+ $, and $\beta = 0.9$. Let $N=2$, $M=10$, $h =1, p= 2$, 
\begin{align*}P =& \begin{bmatrix}   0.4670 &    0.5330 \\    0.4103 &    0.5897
 \end{bmatrix}, 
\\
Q  = &\begin{bmatrix}0.1747  & 0.1716 & 0.1417 & 0.1153 & 0.1095 & 0.0993 & 0.0712 & 0.0658 & 0.0368 & 0.0142 \\
0.0115 & 0.0278 & 0.0537 & 0.0611 & 0.1012 & 0.1176 & 0.1215 & 0.1612 & 0.1667 & 0.1777
 \end{bmatrix}, 
 \end{align*}
and $d = [0\ 1\ 2\ 3 \ 4 \ 8 \ 12 \ 17\ 18\ 19]$.
Then, $\min_{\bm{x}} s^*(\bm{x}) =12$ and $\max_{\bm{x}} s^*(\bm{x}) =17$.
\end{example}
In \autoref{eg4}, \acs{AC} does not hold. \autoref{fig4} presents $v^U_2(\bm{x},s)$ and $v^L_2(\bm{x},s)$ for this example. We note the discontinuity in the expected cost function for two periods, $v^U_2$, obtained by implementing the myopic base stock policy when \acs{AC} does not hold. 

We remark that although $v^U_n(\bm{x},s)$ is piecewise linear in $\bm{x}$ for all $s$ and $n$, in the limit as $n$ approaches $\infty$, we may lose piecewise linearity. Thus, although implementing the policy ``order up to $\max\{s^*(\bm{x}),s\}$'' is straightforward, determining $v^U$, or for that matter $v^U_n$ for large $n$, is computationally demanding. For this reason, we seek an easily computable upper bound on $v^U-v^L$ in the next section. 
\subsection{An Upper Bound on $v^U-v^L$ }
\label{ss:GapUB}
Let $ \Delta = \max_{(\bm{x},s)} \bigg\{ L\big(\bm{x},\max\{s^*(\bm{x}),s\}\big)-L(\bm{x},s^*(\bm{x})) \bigg\} $, and consider the following preliminary result.
\begin{proposition} \label{prop_bd_diff} Assume the policy `order up to $\max\{s^*(\bm{x}),s\}$' is applied and
$min\{s^*(\bm{x}): \bm{x} \in X\} -d_M \leq s(t) \leq \max\{s^*(\bm{x}): \bm{x}\in X\}-d_1$.
Then, for all $\tau \geq t$, 
$$ \min\{s^*(\bm{x}): \bm{x} \in X\} -d_M \leq s(\tau) \leq \max\{s^*(\bm{x}): \bm{x} \in X\} -d_1. $$
\end{proposition}
The next result follows directly from the definition of $\Delta$ and the result above.  
\begin{proposition} \label{prop_delta_def}  Assume the policy `order up to $\max\{s^*(\bm{x}), s\}$' is applied and
$\min\{s^*(\bm{x}) : \bm{x} \in X\} - d_M \leq s(0) \leq \max\{s^*(\bm{x}) : \bm{x} \in X\} - d_1$.  
Then, 
$$ \Delta = \max_{\bm{x}} \{ L \big(\bm{x}, \max\{s^*(\bm{x}): \bm{x} \in X\} - d_1 \big) - L(\bm{x}, s^*(\bm{x}))\}.$$
Thus, assuming $d_M$ is finite, $\Delta$ is finite. 
\end{proposition}

We now present upper bounds on $v_n^U(\bm{x},s)-v^L_n(\bm{x})$ for all $n$ for the special case of backordering systems where $f(y,d) = y-d$ and $c(y,d) = p(d-y)^++h(y-d)^+$. Define $\Delta_n= (1+ \beta+\beta^2+\dots+\beta^{n-1} )\Delta $. Proof of \autoref{difub} follows from a standard induction argument and the fact that the lower bound is independent of $s$. 


\begin{proposition}\label{difub}
Assume $ \min\{s^*(\bm{x}): \bm{x} \in X \} - d_M \leq s(0) \leq \max\{s^*(\bm{x}: \bm{x} \in X \} - d_1$.
For all $\bm{x},s$, and $n$, 
\begin{eqnarray*}  v^U_n(\bm{x},s) - v^L_n(\bm{x}) &\leq& \Delta_n  \
 \text{ and } \ \  v^U(\bm{x},s)-v^L(\bm{x}) \leq \Delta/(1-\beta).
\end{eqnarray*}
\end{proposition}
It follows from \autoref{lemma3}, assuming $s(0) \leq s^*(\bm{x}(0))$, that when \acs{AC} holds, $\Delta = 0$. 

 We remark that $\Delta$ can be determined by a finite set of linear programs. Recall that $\{X_m\}$ is a finite partition of $X$, each element of which is described by a finite set of linear inequalities. Let $d_{\widehat{M}} = \max\{s^*(\bm{x}): \bm{x} \in X\}$, and recall that for $\bm{x} \in X_m$, $s^*(\bm{x}) = d_m$.  Note that $L(\bm{x}, \max\{d_m, d_{\widehat{M}} - d_1\}) - L(\bm{x}, d_m)$ is linear in $\bm{x} \in X_m$ (see \autoref{lemma1}).  Hence, 
\begin{eqnarray} \max\{L(\bm{x}, \max\{d_m, d_{\widehat{M}} - d_1\}) - L(\bm{x}, d_m): \bm{x} \in X_m\} \label{delta_lp}\end{eqnarray}
is a linear program (LP).  When $d_m > d_{\widehat{M}} - d_1$, the optimal objective of this LP is zero (e.g., when $m = \widehat{M}$).  Thus, each element of $\{X_m\}$ either satisfies or does not satisfy the attainability condition for all $\bm{x}$ in that element.  Clearly, $\Delta = 0$ when $\max\{s^*(\bm{x}): \bm{x} \in X\} - d_1 \leq \min\{s^*(\bm{x}): \bm{x} \in X\}$.
The quantity $\Delta$ may be computed as follows by solving the LP defined in \eqref{delta_lp} for each  
 $m \in \mathcal{M}$:
\begin{eqnarray} \Delta = \max_{m \in \mathcal{M}} \max_{\bm{x}}\{L(\bm{x}, \max\{d_m, d_{\widehat{M}} - d_1\}) - L(\bm{x}, d_m): \bm{x} \in X_m\}, \label{eq_comp_del}\end{eqnarray}
where $\mathcal{M} = \{ m:  d_{\widehat{m}}\leq d_m \leq d_{\widehat{M}}-d_1 \} $ and $d_{\widehat{m}} = \min\{s^*(\bm{x}): \bm{x} \in X\}$.
\subsection{Computational Analysis of Myopic Policy for Backordering Systems}
\label{ss:comp_analysis}
We now present a study of the performance of the myopic policy when $\Delta > 0$ and hence \acs{AC} is violated. We consider a set of $216$ test instances of backordering systems with $T=100$ decision epochs that were randomly generated by the procedure described in Appendix \autoref{app_sec4}. Of  the $216$ instances, $\Delta > 0$ for $81$ instances, approximately $37\%$ of the instances. For each instance, on $10,000$ sampled trajectories, we implement the myopic policy to compute a sample average cost $\tilde{v}_T^U(\bm{x}, s)$ and the lower bound policy of resetting the inventory position to the optimal base stock level $s^*(\bm{x})$ of the single period problem (even when doing so requires an infeasible negative reorder quantity) to compute a sample average lower bound cost $\tilde{v}_T^L(\bm{x})$. The algorithms for computing $\tilde{v}_T^U(\bm{x}, s)$ and $\tilde{v}_T^L(\bm{x})$ (Algorithms \ref{algo_ub_tilde} and \ref{algo_lb_tilde}) are presented in Appendix \autoref{algos_policy_impl}. Additionally, we compute $\Delta_T$ over a finite horizon with $T$ decision epochs using the approach outlined in \autoref{ss:GapUB}.

We analyze two quantities: the average observed gap $\tilde{\delta}_T(\bm{x},s) =\tilde{v}_T^U (\bm{x},s) - \tilde{v}^L_T(\bm{x})$ as a fraction of the maximum possible expected difference $\Delta_T$ and as a fraction of the sample average cost of the lower bound policy, where the latter quantity represents a sample average percentage gap approximation.

Table \ref{tab2} presents the results for the $81$ instances in this test set where $\Delta > 0$. The myopic policy generates an average percentage optimality gap of $0.51\%$.  Not unexpectedly, we note that the average observed gap is significantly smaller than the a priori upper bound $\Delta_T$ across this instance set. These results appear to be largely independent of the backorder penalty cost $p$. The average percentage gaps do tend to grow slightly with the number of modulation states N and the number of possible demand outcomes $M$. In summary, however, the quality of the myopic policy, which is optimal when \acs{AC} holds, performs very well on the test instances considered in this numerical study, even when $\Delta>0$. 

\begin{table}[h]
\centering
\begin{tabular}{ccc}
\hline \toprule 
$N$       & $\tilde{\delta}_T(\bm{x},s)/ \Delta_T$ & $\tilde{\delta}_T(\bm{x},s)/ \tilde{v}^L(\bm{x})$ \\
\hline
2       & 3.55\%             & 0.50\%                 \\
3       & 2.15\%             & 0.51\%                 \\
\hline
Overall & 2.65\%             & 0.51\%       		  \\  
\bottomrule  \hline
\end{tabular}
\hspace{4ex}
\begin{tabular}{ccc}
\hline \toprule 
$M$       & $\tilde{\delta}_T(\bm{x},s)/ \Delta_T$ & $\tilde{\delta}_T(\bm{x},s)/ \tilde{v}^L(\bm{x})$ \\
\hline
3       & 1.68\%             & 0.41\%                 \\
4       & 3.06\%             & 0.45\%                 \\
5       & 3.02\%             & 0.60\%                 \\
\bottomrule  \hline
\end{tabular}%
\vspace{4ex}
\begin{tabular}{ccc}
\hline \toprule 
$p$       & $\tilde{\delta}_T(\bm{x},s)/ \Delta_T$ & $\tilde{\delta}_T(\bm{x},s)/ \tilde{v}^L(\bm{x})$ \\
\hline
1.5     & 2.42\%             & 0.58\%                 \\
2       & 2.43\%             & 0.45\%                 \\
3       & 3.23\%             & 0.48\%                 \\
\bottomrule  \hline
\end{tabular}

\caption{Performance of the myopic policy on subset of test instances for which \acs{AC} is violated. \label{tab2}} 
\end{table}

\section{Reorder Cost Case \label{sec5}}
We now consider the case where there is a reorder cost $K\geq 0$.
The following results combine the ideas presented for the $K=0$ case with straightforward extensions of earlier results in the literature. Let the operator $H^K$ be defined as:
$$ [H^Kv](\bm{x},s) = \min_{y\geq s} \bigg\{ K\xi(y-s) + [Gv](\bm{x},y) \bigg\},$$
where $\xi(k) = 0$ if $k=0$ and $\xi(k) = 1$ if $k\neq0$ and 
$$[Gv](\bm{x},y) = L(\bm{x},y) + \beta \sum_{d,z} \sigma(d,z,\bm{x})v\big(\bm{\lambda}(d,z,\bm{x}),f\big(y, d\big)\big). $$
We note that when $K=0$, $H^K = H$, as defined in \autoref{model}. 

We now assume that $K>0$. Our objective is to present conditions under which $(s,S)$ policies exist and how such policies can be computed. 
\subsection{$K$-convexity and Optimal $(s, S)$ Policies} 
We now present our first result following a key definition: 
the real-valued function $g$ is $K$-convex if for any $s\leq s'$, 
$$ g\big(\lambda s + (1-\lambda)s'\big) \leq \lambda g(s) + (1-\lambda) \big( g(s') + K \big), \text{ for all } \lambda \in [0, 1].$$
Proof of the following result is a direct extension of results in \cite{scarf60} and elsewhere. 
\begin{proposition}\label{prop6}
Assume: (i) $v(\bm{x},s)$ is $K$-convex in $s$ for all $\bm{x}$, (ii) $f(y,d)$ is non-decreasing and convex in $y$ for all $d$, and (iii) $c(y,d)$ is convex in $y$ and $\lim_{\abs{y} \to \infty} c(y,d) \to \infty $ for all $d$. 
Then,
\begin{description}
\item[1.] $[Gv](\bm{x},y)$ is $K$-convex in $y$ for all $\bm{x}$,
\item[2.] $[H^Kv](\bm{x},s)$ is $K$-convex in $s$ for all $\bm{x}$, and
\item[3.] $[H^Kv](\bm{x},s) = \begin{cases}K+[Gv]\big(\bm{x},S^*(\bm{x},v)\big) & s\leq s^*(\bm{x}, v) \\ [Gv](\bm{x}, s) & \text{otherwise,} \end{cases}$
\\where:
 $S^*(\bm{x}, v)$ is the smallest integer minimizing $[Gv](\bm{x}, y)$ with respect to $y$, and $s^*(\bm{x},v)$ is the smallest integer such that $[Gv]\big(\bm{x},s^*(\bm{x},v)\big) \leq K+[Gv]\big(\bm{x}, S^*(\bm{x},v)\big).$
\end{description}
\end{proposition}
Thus, the fact that $v(\bm{x},s)$ is $K$-convex and non-decreasing in $s$ for all $\bm{x}$ leads to the existence of an optimal policy that is of $(s, S)$ form: if the inventory drops below $s$, then order up to $S$; otherwise, do not replenish.
\subsection{Bounds on {$s_n$ and $S_n$}}
Let $v_0=0$, $v_{n+1} = H^Kv_n$ for all $n\geq 0$, and $G_n(\bm{x},y) = [Gv_n](\bm{x}, y)$. Let $S_n(\bm{x})$ be the smallest integer such that $G_n \big(\bm{x}, S_n(\bm{x})\big)\leq G_n(\bm{x}, y)$ for all $y$, and let $s_n(\bm{x})$ be the smallest integer such that $G_n\big(\bm{x}, s_n(\bm{x})\big) \leq K+G_n\big(\bm{x}, S_n(\bm{x})\big)$. Following \cite{veinott65c}, we now define four real-valued functions that represent bounds on the set $\big\{ \big(s_n(\bm{x}),S_n(\bm{x})\big): n\geq 0  \big\}$. Let the values $\underline{s}(\bm{x}), \overline{s}(\bm{x}), \underline{S}(\bm{x})$, and $\overline{S}(\bm{x})$ be the smallest integers such that: 
\begin{eqnarray}
L\big(\bm{x},\underline{S}(\bm{x})\big) &\leq& L(\bm{x},y) \ \forall \ y
\\L\big(\bm{x}, \underline{s}(\bm{x}) \big) &\leq & K+L\big(\bm{x},\underline{S}(\bm{x})\big)
\\\beta K + L\big(\bm{x}, \underline{S}(\bm{x})\big) &\leq& L\big(\bm{x},\overline{S}(\bm{x}) \big), \ \overline{S}(\bm{x}) \geq \underline{S}(\bm{x})
\\ L\big(\bm{x},\overline{s}(\bm{x})\big) &\leq& L\big(x,\underline{S}(\bm{x})\big)+(1-\beta)K,
\end{eqnarray}
where, from earlier results, $\underline{S}(\bm{x})$ can be restricted to the set $\{d_1, \dots, d_M\}$ and where $\underline{S}$ is identical to the functions $s^*$ and $S_0$. We remark that the convexity of $L(\bm{x},y)$ in $y$ for all $\bm{x}$ insures that for each $\bm{x}$, $\underline{s}(\bm{x}) \leq \overline{s}(\bm{x})\leq \underline{S}(\bm{x}) \leq \overline{S}(\bm{x})$.
\subsection{A Partition based on {$(\underline{s}, \overline{s}, \underline{S}, \overline{S})$}}
 Extending results in \cite{veinott65c}, we now show that for all $\bm{x}$ and $n$, $\underline{s}(\bm{x}) \leq s_n(\bm{x}) \leq\overline{s}(\bm{x})$ and  $\underline{S}(\bm{x}) \leq S_n(\bm{x}) \leq\overline{S}(\bm{x})$  and that for the  infinite horizon discounted case,  $\underline{s}(\bm{x}) \leq s^*(\bm{x}) \leq\overline{s}(\bm{x})$ and  $\underline{S}(\bm{x}) \leq S^*(\bm{x}) \leq\overline{S}(\bm{x})$, where $(s^*, S^*)$ represents an $(s, S)$ belief-dependent optimal policy. Proof is presented in the e-companion of this paper. 
 
 \begin{proposition} \label{prop7}
\begin{description}
\item[(a)] For the $n$-period problem, for all $\bm{x}$, there exists an optimal $(s, S)$ policy $(s_n(\bm{x})$, $S_n(\bm{x}))$,  where $\underline{s}(\bm{x}) \leq s_n(\bm{x}) \leq \overline{s}(\bm{x}) \leq \underline{S}(\bm{x}) \leq S_n(\bm{x}) \leq \overline{S}(\bm{x})$. 
\item[(b)] For the infinite horizon problem, for all $\bm{x}$ there is an epoch-invariant $(s,S)$ policy $(s^*(\bm{x}), S^*(\bm{x}))$, where $\underline{s}(\bm{x}) \leq s^*(\bm{x}) \leq \overline{s}(\bm{x}) \leq \underline{S}(\bm{x}) \leq S^*(\bm{x}) \leq \overline{S}(\bm{x})$. 
\end{description}
\end{proposition}

Assume $f(y,d) = y-d$ and $c(y,d) = p(d-y)^++h(y-d)^+ $ and recall from \autoref{lemma2} that $S^*(\bm{x}) = d_m$ if $\bm{x}$ satisfies 
\begin{eqnarray} \label{eq5.1}
 \bm{x}\sum_{k=1}^{m-1}\bm{P}(d_k)\underline{1} < \frac{p}{p+h} \leq \bm{x} \sum_{k=1}^m \bm{P}(d_k)\underline{1}. 
\end{eqnarray}
Given $\underline{S}(\bm{x}) = d_m$, let $\underline{s}(\bm{x}) = d_i$, $\overline{s}(\bm{x}) = d_j$, and $\overline{S}(\bm{x}) = d_n$ satisfy 
\begin{eqnarray}
A_i(\bm{x})d_i+B_i(\bm{x}) &\leq& K+A_m(\bm{x})d_m + B_m(\bm{x})
\\ A_j(\bm{x}) d_j + B_j(\bm{x}) &\leq& (1-\beta)K+A_m(\bm{x})d_m+B_m(\bm{x})
\\ \label{eq5.4} \beta K+ A_m(\bm{x})d_m + B_m(\bm{x}) &\leq& A_n(\bm{x}) d_n + B_n(\bm{x}).
\end{eqnarray}
Let $\overline{X}(\underline{s}, \overline{s}, \underline{S}, \overline{S})$ be the set of all $\bm{x}\in {X}$ such that $\underline{s} = d_i$, $\overline{s} = d_j$, $\underline{S} = d_m$, and $\overline{S} = d_n$ are the smallest integers satisfying the five linear inequalities in Eqs. \ref{eq5.1} - \ref{eq5.4}. Note that the set of all $\overline{X}(\underline{s}, \overline{s}, \underline{S}, \overline{S})$ such that $\overline{X}(\underline{s}, \overline{s}, \underline{S}, \overline{S})$ is non-null is a partition of ${X}$.  
\begin{example}\label{eg7}
Consider \autoref{eg1} with reorder cost, $K=5$. 
\begin{figure}[h]
\centering
\includegraphics[scale = 0.6, trim={65pt 65pt 0pt 0pt},clip]{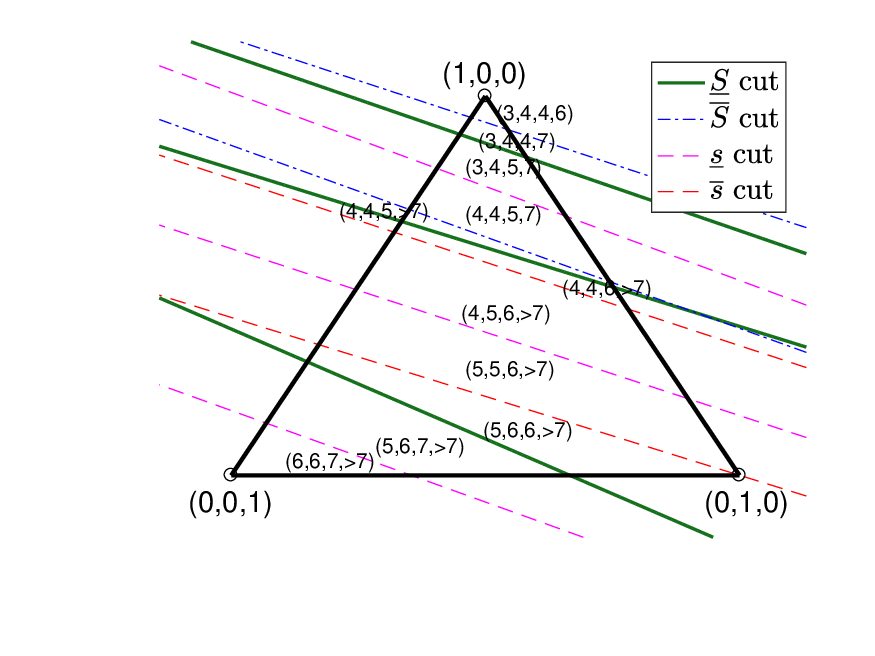}
\caption{Example of partition of ${X}$, with $N=3$ and $M=7$} \label{fig6}
\end{figure}
Each region in the triangle is described by $(i, j, m, n)$ where, $(\underline{s}, \overline{s}, \underline{S}, \overline{S}) = \big(d(i), d(j), d(m), d(n)\big)$. For example, the region labeled as $(4,4,5,>7)$ in \autoref{fig6} has $ (\underline{s}, \overline{s}, \underline{S}, \overline{S}) = (20, 20, 25, 36)$. This implies that $s^*(\bm{x}) =d_4= 20, \ \forall \bm{x}\in (4,4,5,>7)$. The search interval for $S^*(\bm{x})$ is also significantly restricted to the demand outcomes between $\underline{S}$ and $\overline{S}$, making the computation very easy. We remark that it is possible $\overline{S} > d_M$, as indicated (by $>7$) in \autoref{fig6}. The corresponding $\overline{S}$ is $36$ in $X_5$ and $X_6$ and it equals $38$ in $X_7$.
\end{example}  

A description of the determination of the sets $\Gamma_n(s)$, where $v_n(\bm{x},s) = \min \{ \bm{x}\bm{\gamma}: \bm{\gamma} \in \Gamma_n(s)\}$, can be found in the e-companion of this paper. 


\section{Conclusions\label{sec6}} 
We have presented and analyzed an inventory control problem having a modulation process that affects demand and that is partially observed by the demand and AOD processes.  Assuming  \acs{AC} holds and the reorder cost $K = 0$, we generalized results found throughout the literature that there exists an optimal policy that is a myopic base stock policy. We also developed a simple, easily implemented description of the optimal myopic base stock levels, as a function of the belief function.  
When \acs{AC} is violated and $K=0$, we examined the the base stock policy that is optimal when \acs{AC} holds as a suboptimal policy yielding an upper bound. We presented a lower bound on the optimal expected cost function and a bound on the difference between the upper and lower bounds. A numerical study indicated that the bound on the difference between these two bounds can be quite small, indicating that even when \acs{AC} is violated, the optimal base stock policy for the case where \acs{AC} is not violated may be quite good.  

When $K > 0$, we showed that there exists optimal $(s, S)$ policies, dependent on the belief function, and determined upper and lower bounds on $s$ and $S$ for the finite and infinite horizon cases, where each bound is dependent on the belief function of the modulation process. We showed that each of these bounds and the values of $s$ and $S$ for the finite and infinite horizon cases are constant within regions of the belief space and that these regions can be described by a finite number of linear inequalities. We outlined an approach for determining an optimal $(s, S)$ policy and the resultant expected cost function for the finite horizon case. 

Future work may focus on extending the approach of formulation and solution presented for the current problem to other problems with similar structure. Another direction of potential interest would be the development of heuristic approaches for the case when \acs{AC} does not hold, as pursued by \citet{malladi21} for a distributed production system with mobile production capacity.  

\bibliographystyle{spbasic-no-url}      
\bibliography{pomdp_inventory} 

\newpage
\section*{Appendix}
\renewcommand{\thesubsection}{A\arabic{subsection}} 
\subsection{Proof of Result in Section 2}
\begin{proof}[Proof of \autoref{lemma2}.]
If $s^*(\bm{x}) = d_m$, then
\begin{eqnarray*}
A_{m-1}(\bm{x})d_{m-1} + B_{m-1}(\bm{x}) &>& A_m(\bm{x})d_m+ B_m(\bm{x}),
\\ A_{m+1}(\bm{x})d_{m+1} + B_{m+1}(\bm{x}) &\geq& A_m(\bm{x})d_m+B_m(\bm{x}),
\end{eqnarray*}
which leads to the result.  
\end{proof}
\subsection{Proofs of Results in Section 3 \label{app_sec3}}
\begin{proof}[Proof of \autoref{prop1}.]
By induction. Letting $v_0 = 0$, we note that 
$$
v_1(s,\bm{x})=\min_{y\geq s} L(\bm{x},y) = L\big(\bm{x},\max\{s^*(\bm{x}),s\}\big)
$$
for all $\bm{x}$ and $L\big(\bm{x},\max\{s^*(\bm{x}),s\}\big)$ is non-decreasing and convex in $s$. Thus, the result holds true for $n=1$ (and, trivially for $n=0$). Assume the result holds for $n$. Then, for $s\leq s^*(\bm{x})$,
\begin{eqnarray*}
&&v_{n+1}(\bm{x},s) \leq L\big(\bm{x},s^*(\bm{x})\big) + \beta \sum_{d,z} \sigma(d,z,\bm{x})v_n\big(\bm{\lambda}(d,z,\bm{x}) , f\big( s^*(\bm{x}), d\big)\big)
\\ &&= L\big(\bm{x},s^*(\bm{x})\big) + \beta \sum_{d,z} \sigma(d,z,\bm{x})v_n\big(\bm{\lambda}(d,z,\bm{x}), s^*\big(\bm{\lambda}(d,z,\bm{x})\big) \big) \ \text{(using \autoref{A1}).}
\end{eqnarray*}
\begin{eqnarray*}
\text{Also, } v_{n+1}(\bm{x},s) &\geq & \min_{y\geq s} L(\bm{x},y)+ \beta \sum_{d,z} \sigma(d,z,\bm{x}) \min_y v_n \big (\bm{\lambda}(d,z,\bm{x}),f\big(y, d\big) \big)
\\ &=& L\big(\bm{x},s^*(\bm{x})\big) +\beta \sum_{d,z} \sigma(d,z,\bm{x}) v_n \big (\bm{\lambda}(d,z,\bm{x}), s^*\big(\bm{\lambda}(d,z,\bm{x})\big) \big)
\\&=& L\big(\bm{x},s^*(\bm{x})\big) +\beta \sum_{d,z} \sigma(d,z,\bm{x})  v_n \big (\bm{\lambda}(d,z,\bm{x}), f\big( s^*(\bm{x}), d\big) \big).
\end{eqnarray*}
$\text{ Thus, for } s\leq s^*(\bm{x}), \text{ } $
$$ v_{n+1}(\bm{x},s) = L\big(\bm{x},s^*(\bm{x})\big) + \beta \sum_{d,z} \sigma(d,z,\bm{x}) v_{n}\big(\bm{\lambda}(d,z,\bm{x}), f \big(s^*(\bm{x}),d\big)\big), $$
and $ v_{n+1}(\bm{x},s) = v_{n+1}(\bm{x}, s^*(\bm{x})) .$

Assume $s \geq s^*(\bm{x})$. Note 
$$ v_{n+1}(\bm{x},s) \leq L(\bm{x},s) + \beta \sum_{d,z} \sigma(d,z,\bm{x}) v_n\big(\bm{\lambda}(d,z,\bm{x}),f\big(s, d\big)\big). $$
\begin{eqnarray*}
\text{ Also, } v_{n+1}(\bm{x},s) &\geq& \min_{y\geq s} L(\bm{x},y) + \beta \sum_{d,z} \sigma(d,z,\bm{x})\min_{y\geq s} v_n \big(\bm{\lambda}(d,z,\bm{x}),f\big(y, d\big)\big)
\\ &=& L(\bm{x},s) + \beta \sum_{d,z} \sigma(d,z,\bm{x}) v_n\big(\bm{\lambda}(d,z,\bm{x}), f \big(s, d \big)\big),
\end{eqnarray*}
$$\text{ and hence for } s\geq s^*(\bm{x}), \ \text{ } v_{n+1}(\bm{x},s) = L(\bm{x},s) + \beta\sum_{d,z} \sigma(d,z,\bm{x}) v_n\big(\bm{\lambda}(d,z,\bm{x}),f \big(s, d \big)\big) $$
and is non-decreasing and convex in $s$. 
\end{proof}

\begin{proof}[{Proof of \autoref{lemma6}.}]
It is sufficient to show that if $y\leq y'$ and $\bm{x} \preceq \bm{x'}$, then,
$$L(\bm{x},y) - L(\bm{x},y') \leq L(\bm{x'},y)-L(\bm{x'},y'),$$
which follows from the assumptions and \cite[Lemma 4.7.2]{puterman94}.
\end{proof}

Ideally, we would want to select $\widehat{\bm{x}}^{d,z}$ so that $s^*(\bm{x'})\leq s^*(\widehat{\bm{x}}^{d,z})$ for all $\bm{x'}$ such that $\bm{x'} \preceq \bm{\lambda}(d,z,\bm{x}) \ \ \forall \ \bm{x}\in X$, for all $(d,z)$, which would strengthen \autoref{lemma6} as much as possible. We construct such an $\widehat{\bm{x}}^{d,z}$ after the following preliminary result.  

\begin{lemma} \label{lemma7a}
The set $\{ \bm{\lambda}(d,z,\bm{x}): \bm{x}\in X\}$ = $\bigg\{ \sum_i \xi_i \bm{\lambda}(d,z,\bm{e_i}) : \xi_i \geq 0 \ \forall i, \sum_i \xi_i = 1 \bigg \}.$
\end{lemma}
We remark that if $\bm{x} \preceq \bm{x'}$ and $\bm{x} \preceq \bm{x''}$, then $\bm{x}\preceq\alpha \bm{x'}+(1-\alpha)\bm{x''} $ for all $\alpha \in [0,1]$. Thus, if $\widehat{\bm{x}}^{d,z}$ is such that $\widehat{\bm{x}}^{d,z} \preceq \lambda(d,z,\bm{e_i})$ for all $i$, then $\widehat{\bm{x}}^{d,z}$ is such that $\widehat{\bm{x}}^{d,z} \preceq \bm{x'}$ for all $\bm{x'} \in \big \{ \bm{\lambda}(d,z,\bm{x}) : \bm{x}\in X \big\}$.  

\subsection*{Construction of $\widehat{x}^{d,z}$}
We now construct  $\widehat{\bm{x}}^{d,z}$. 
Let 
\begin{eqnarray*}
\widehat{x}_N^{d,z} &=& \min \big\{ \bm{\lambda}_N(d,z,\bm{e_i}), i=1,\dots,N \big\}
\\ \widehat{x}_n^{d,z} &=& \min \bigg \{ \sum_{k=n}^N \bm{\lambda}_k(d,z,\bm{e_i}), i =1,\dots, N \bigg\} -\sum_{k=n+1}^N \widehat{x}_k^{d,z}, \ \  n= N-1, \dots, 2
\\ \widehat{x}_1^{d,z} &=& 1 - \sum_{k=2}^N \widehat{x}_k^{d,z}.
\end{eqnarray*}
By construction, $\widehat{\bm{x}}^{d,z}\preceq \bm{\lambda}(d,z, \bm{x}) \ \forall \ \bm{x}\in X$. We now show that $\widehat{\bm{x}}^{d,z}\in X$ and that $s^*(\bm{x'})\leq s^*(\widehat{\bm{x}}^{d,z})$ for all $\bm{x'}\in X$ such that $\bm{x'}\preceq \bm{\lambda}(d,z,\bm{x}) \ \forall \ \bm{x}\in X$. 

\begin{lemma} \label{lemma8} (i) $\widehat{\bm{x}}^{d,z} \in X$. (ii) For any $\bm{x'}\preceq \bm{\lambda}(d,z,\bm{x}) \ \forall \ \bm{x}\in X, s^*(\bm{x'}) \leq s^*(\widehat{\bm{x}}^{d,z})$. 
\end{lemma}

\begin{proof}[{Proof of \autoref{lemma8}.}]
We have the following:
\begin{description}
\item[(i)] Clearly, $0\leq \widehat{\bm{x}}^{d,z}_N \leq 1$ and $\sum_{n=1}^N \widehat{\bm{x}}_n^{d,z} = 1$. It is sufficient to show $0\leq \widehat{\bm{x}}^{d,z}_n, n=N-1, \dots, 1$. Note 
$$\sum_{k=n+1}^N\widehat{\bm{x}}^{d,z}_k = \min_{1\leq i\leq N} \bigg \{\sum_{k=n+1}^N \lambda_k(d,z,\bm{e_i})\bigg\} \leq \sum_{k=n+1}^N \lambda_k(d,z,\bm{e_i}) \leq \sum_{k=n}^N \lambda_k(d,z,\bm{e_i}), \ \forall \ i.$$
Thus, $\sum_{k=n+1}^N\widehat{x}^{d,z}_k \leq \min_{1\leq i\leq N} \bigg \{\sum_{k=n}^N \lambda_k(d,z,\bm{e_i})\bigg\} = \sum_{k=n}^N\widehat{x}^{d,z}_k$, 
and hence $\widehat{x}_n^{d,z}\geq 0$. 
\item[(ii)] Let $\bm{x'}\preceq \bm{\lambda}(d,z,\bm{x}) \ \forall \ \bm{x} \in X$ and assume $s^*(\widehat{\bm{x}}^{d,z}) < s^*(\bm{x'})$. Then, there is an $n\in \{1, \dots, N\}$ such that $\sum_{k=n}^N x_k' > \sum_{k=n}^N \widehat{x}_k^{d,z}$. However, 
$\sum_{k=n}^N \widehat{x}^{d,z}_k = \min_{1\leq i\leq N}$ $\bigg \{\sum_{k=n}^N \lambda_k(d,z,\bm{e_i})\bigg\}$, which leads to a contradiction of the assumption that $\bm{x'} \preceq \bm{\lambda}(d,z,\bm{x})$ $\forall \bm{x} \in X$. 
\end{description}
\end{proof}

\subsection*{\label{sec3.3}Computing the Expected Cost Function, $v_n$}
We now present a procedure for computing $v_n(s, \bm{x})$. We only consider the case where $s=s^*(\bm{x})$ due to \autoref{prop1} and \autoref{lemma3}. For notational simplicity, we assume that $\text{Pr}\big(z(t+1) \mid \mu(t+1),\mu(t)\big)$ is independent of $\mu(t+1)$ and $\mu(t)$. Extension to the more general case is straightforward.

Assume $v_0=0$, let $n=1$, and recall $v_1\big(\bm{x},s^*(\bm{x})\big) = L\big(\bm{x},s^*(\bm{x})\big). $
Note  $ L(\bm{x},y) = \bm{x}\overline{\bm{\gamma}}_y, \ \text{ where } \overline{\bm{\gamma}}_y = \sum_{d,z} \bm{P}(d,z)\ \underline{1}\ c(y,d)$. Let $\Gamma_1 = \{ \overline{\bm{\gamma}}_y\}$, and note that if $c(y, d) = p(d-y)^+ + h(y-d)^+$, it is sufficient to consider only $y \in \{d_1,\dots, d_M\}$.  Then,
 $v_1\big(\bm{x},s^*(\bm{x})\big)= \min \big\{\bm{x}\overline{\bm{\gamma}}:\overline{\bm{\gamma}} \in \Gamma_1 \big\}.$  Assume there is a finite set $\Gamma_n$ such that  $ v_n \big(\bm{x}, s^*(\bm{x})\big) = \min \big\{\bm{x}\bm{\gamma}: \bm{\gamma}\in \Gamma_n \big\}.$ Then,
\begin{eqnarray*}
v_{n+1}\big(\bm{x},s^*(\bm{x}) \big) &=& L\big(\bm{x},s^*(\bm{x})\big) + \beta \sum_{m=1}^M \sigma(d_m,\bm{x}) v_n\big(\bm{\lambda}(d_m, \bm{x}),f\big(s^*(\bm{x}),d_m\big)\big)
\\ &=&\min \big\{\bm{x}\overline{\bm{\gamma}}:\overline{\bm{\gamma}} \in \Gamma_1 \big\} + \beta \sum_{m=1}^M \sigma(d_m,\bm{x}) v_n\big(\bm{\lambda}(d_m,\bm{x}),s^*(\bm{\lambda}(d_m,\bm{x}))\big)
\\&=& \min \big\{\bm{x}\overline{\bm{\gamma}}:\overline{\bm{\gamma}} \in \Gamma_1 \big\} + \beta \sum_{m=1}^M \sigma(d_m,\bm{x})  \min \big\{\bm{\lambda}(d_m,\bm{x})\bm{\gamma}: \bm{\gamma}\in \Gamma_n \big\}
\\&=& \min_{\overline{\bm{\gamma}}} \min_{\bm{\gamma_1}}\dots \min_{\bm{\gamma_M}} \bigg\{ \bm{x}\overline{\bm{\gamma}} + \beta \sum_{m=1}^M \sigma(d_m,\bm{x})\bm{\lambda}(d_m,\bm{x}) \bm{\gamma_m} \bigg \}
\\&=& \min_{\overline{\bm{\gamma}}} \min_{\bm{\gamma_1}}\dots \min_{\bm{\gamma_M}} \bigg\{ \bm{x} \bigg[\overline{\bm{\gamma}} + \beta \sum_{m=1}^M \bm{P}(d_m) \bm{\gamma_m} \bigg] \bigg\}
\end{eqnarray*}
Thus, $\Gamma_{n+1}$ is the set of all $\bm{\gamma}$ such that $ \bm{\gamma} = \overline{\bm{\gamma}} + \beta \sum_{m=1}^M \bm{P}(d_m)\bm{\gamma_m}, $
where $\overline{\bm{\gamma}} \in \Gamma_1$ and $\bm{\gamma_m} \in \Gamma_n$ for all $m =1,\dots, M$, and for all $n$, $v_n\big(\bm{x}, s^*(\bm{x})\big) $ is piecewise linear and concave in $\bm{x}$. 

Let $\abs{\Gamma}$ be the cardinality of the set $\Gamma$. Then, $\abs{\Gamma_{n+1}} = \abs{\Gamma_1} \abs{\Gamma_n}^M$, where $\abs{\Gamma_1} \leq M$, and hence the cardinality of $\Gamma_n$ can grow rapidly. Many of the vectors in the sets $\Gamma_n$ are redundant and can be eliminated, reducing both computational and storage burdens. An exhaustive literature study of elimination procedures and other solution methods for solving POMDPs can be found in \cite{chang15a}.

 \subsection{Proofs of Results in Section 4 \label{app_sec4}}
\begin{proof}[{Proof of \autoref{lemma8.5}}]
Assume $f(y,d) = y-d$ and $c(y,d) = p(d-y)^++h(y-d)^+ $, recall that elements of $\mathcal{P}_1$ are sets of the form $\{\bm{x}\in X: s^*(\bm{x}) = d_m\}$ for all $d_m$ such that $\min_{\bm{x}} s^*(\bm{x}) \leq d_m \leq \max_{\bm{x}} s^*(\bm{x})$. Further recall that $\{\bm{x}\in X: s^*(\bm{x}) = d_m\}$ is the set of all $\bm{x}$ such that 
$$\sum_{k=1}^{m-1} \sigma(d_k,\bm{x}) <  p/(p+h) \leq \sum_{k=1}^m \sigma(d_k,\bm{x}), $$
\text{ or equivalently, } $$\quad \quad \quad \quad \bm{x}\sum_{k=1}^{m-1}\bm{P}(d_k)\underline{1} <  p/(p+h) \leq \bm{x}\sum_{k=1}^m \bm{P}(d_k)\underline{1},$$ 
 which represents two linear inequalities. Further, for $\bm{x}\in \{ \bm{x} \in X: s^*(\bm{x}) =d_m\}$,
$ v_1^U(\bm{x},s) = A_l(\bm{x})d_l+B_l(\bm{x}) $ 
for $l=\max\{s^*(\bm{x}),s\}$, where we note 
$$  A_l(\bm{x})d_l+B_l(\bm{x}) = \bm{x} \big [h\sum_{k=1}^l (d_l-d_k)\bm{P}(d_k)\underline{1} + p \sum_{k=l+1}^M(d_k-d_l)\bm{P}(d_k)\underline{1} \big], $$
where $A_j(\bm{x})$ and $B_j(\bm{x})$ are defined in \autoref{lxy}. Thus, on each element of $\mathcal{P}_1$, $v_1^U$ is linear in $\bm{x}$ for each $s$ and each element of $\mathcal{P}_1$ is described by a finite number of linear inequalities. 

Let $(\bm{x},s)$ be such that $d_l \leq \max\{s^*(\bm{x}),s\} \leq d_{l+1}$ for all $\bm{x}$ in an element $\{\bm{x}\in X: s^*(\bm{x})=d_m\}$. Further, let $d_{l(d,z)} \leq  \max\{s^*(\bm{\lambda}(d,z,\bm{x})),\max\{s^*(\bm{x}),s\}-d\} \leq d_{l(d,z) + 1}$ for all $\bm{x}$ in an element $\{\bm{x}\in X: s^*(\bm{\lambda}(d,z,\bm{x})) = d_{m(d)}\}$, which is the set of all $\bm{x}$ such that:
\begin{eqnarray*}
\bm{\lambda}(d,z,\bm{x})\sum_{k=1}^{m(d)-1}\bm{P}(d_k)\underline{1} < p/(p+h) \leq \bm{\lambda}(d,z,\bm{x})\sum_{k=1}^{m(d)} \bm{P}(d_k)\underline{1},
\end{eqnarray*}
or equivalently, for all $\bm{x}$ such that  $\sigma(d,\bm{x}) \neq 0$, 
\begin{eqnarray*}
\bm{x}\bm{P}(d,z)\sum_{k=1}^{m(d)-1}\bm{P}(d_k)\underline{1} < \big(p/(p+h)\big)\bm{x}\bm{P}(d,z)\underline{1} \leq \bm{x}\bm{P}(d,z)\sum_{k=1}^{m(d)} \bm{P}(d_k)\underline{1},
\end{eqnarray*}
where we assume $m$ and $m(d)$ for all $d$ have been chosen so that the finite set of linear inequalities describes a non-null subset of $X$. We note that for such a subset, 
\begin{eqnarray*}
v^U_{n+1}(\bm{x},s) &=& A_l(\bm{x})d_l +B_l(\bm{x})
\\ &&+ \beta \sum_{d,z} \sigma(d,z,\bm{x}) \bigg[A_{l(d,z)}(\bm{\lambda}(d,z,\bm{x}))d_{l(d,z)} + B_{l(d,z)}(\bm{\lambda}(d,z,\bm{x})) \bigg]
\\ &=& \bm{x}\bigg[h\sum_{k=1}^l(d_l-d_k)\bm{P}(d_k)\underline{1} + p \sum_{k=l+1}^M(d_k-d_l)\bm{P}(d_k)\underline{1}
\\&&+\beta\sum_d \bigg[h \sum_z\sum_{k=1}^{l(d,z)}(d_{l(d,z)}-d_k)\bm{P}(d,z)\bm{P}(d_k)\underline{1} 
\\&&+ p \sum_z \sum_{k=l(d,z)+1}^N(d_k-d_{l(d,z)})\bm{P}(d,z)\bm{P}(d_k)\underline{1}\bigg]\bigg].
\end{eqnarray*}
The resulting partition $\mathcal{P}_2$ is at least as fine as $\mathcal{P}_1$ and each element in $\mathcal{P}_2$ is described by a finite set of linear inequalities.
 We have shown that on each element in $\mathcal{P}_2$, $v_2^U(\bm{x},s)$ is linear in $\bm{x}$ for each $s$. A straightforward induction argument shows these characteristics hold for all $n$. We illustrate by example (through \autoref{eg4})  how $v_n^U(\bm{x},s)$ may be discontinuous in $\bm{x}$ for fixed $s$. 
 \end{proof}
 \begin{proof}[{Proof of \autoref{prop_bd_diff}.}] It is sufficient to show the result holds for $\tau = t+1$.  There are two cases.  First, let $s(t) \leq s^*(\bm{x}(t))$.  Then, $s(t+1) = s^*(\bm{x}(t))-d(t)$.  We note
\begin{eqnarray*}
\min\{s^*(\bm{x}): \bm{x} \in X\}- d_M &\leq& s^*(\bm{x}(t))- d(t) 
\\ &\leq& \max\{s^*(\bm{x}): \bm{x} \in X\}- d_1 \text{ and hence,} 
\\ \min\{s^*(\bm{x}): \bm{x} \in X\} -d_M &\leq& s(t+1) \leq \max\{s^*(\bm{x}): \bm{x} \in X\}- d_1
\end{eqnarray*}
Second, let $s^*(\bm{x}(t)) \leq s(t)$.  Then, $s(t+1) = s(t)- d(t) \geq s^*(\bm{x}(t))-d(t)$.  We note 
\begin{eqnarray*}
\min\{s^*(\bm{x}): \bm{x} \in X\} - d_M \leq s^*(\bm{x}(t)) - d(t) \leq s(t) - d(t) 
\\ \leq \max\{s^*(\bm{x}): \bm{x} \in X\} - d_1, \text{ and hence, } 
\\ \min\{s^*(\bm{x}): \bm{x} \in X\} -d_M \leq s(t+1) \leq \max\{s^*(\bm{x}): \bm{x} \in X\} - d_1.
\end{eqnarray*} 
\end{proof}

\subsection{Design of Instances for Computational Study}
\label{appendix:instances}
We describe the generation of computational instances for Section \autoref{ss:comp_analysis}. Each instance describes a backordering system with no fixed ordering cost. For each combination of number of modulation states $N \in \{2,3\}$, number of demand outcomes $M \in \{3,4,5\}$, randomly generate $M$ unique ordered integer demand outcomes from $[0, D_L]$ for each $D_L \in \{20, 100, 250, 500, 750, 1000\}$. Set the lowest demand outcome $d(0) = 0$ ( to encourage A1 violation), randomly sample probability transition matrix $\{P(i,j)\}$ and probability mass function for each modulation state  $\{Q(d,j)\}$ such that the $N$ ordered expected demands $ED_i, i =1, \dots, N$ are quite distinct and satisfy: 
\begin{itemize}\item	$ED_1<= 0.5 d(M)$ and $ED_2 > 0.5 d(M)$ and $ED_2-ED_1 > 0.25 d(M)$ OR $ED_2 - ED_1 > 0.5 d(M)$, when $N = 2$ and \item	$ED_1 <= 0.4 d(M)$ and $ED_2 > 0.4 d(M)$ and $ED_2 <= 0.7 d(M)$ and $ED_3 > 0.7 d(M)$ and $ED_2 - ED_1 > 0.2 d(M)$ and $ED_3 – ED_2 > 0.2 d(M)$, when $N = 3$.\end{itemize}
Set the number of decision epochs $T$ to $100$ and vary backorder cost per unit per period $p$ as $\{1.5, 2, 3\}$, while keeping the holding cost $h$ at $1$.

 \subsection{Algorithms for Computational Study \label{algos_policy_impl}}
 \begin{minipage}{\linewidth}
\begin{algorithm}[H]
\begin{algorithmic}[1]
 \caption{Sample Average Cost of Myopic Policy $\tilde{v}^U_T(\bm{x}, s) $} \label{algo_ub_tilde}
 \Procedure{$\tilde{v}^U_T$}{$\bm{x},s$}
 \For{$i = 1,\dots,N $} 
 \State Set trajectory counter $n \gets 1$. 
 \While {$n<nSamples$ }
 \State Initialize $(\bm{x}, s) \gets (\bm{e}_i, 0)$. Initialize partially observed modulation state ${\mu_0} \gets \mu^i$.
 \For {$t=1,\dots,T$} 
\State Calculate $y = \max\big( s, s^*(x)  \big)$. \State Randomly sample the current modulation state $\mu$ with probability $\bm{e}(\mu_0) P_{ij}$. 
 \State Randomly sample the current demand outcome $d^n_t$ with probability $\text{Pr}(d \mid \mu)$.
 \State Compute $c_t^n = h(y-d^n_t)^+ + p(d^n_t-y)^+$.
 \State Update $\mu_0 \gets \mu $, $\bm{x} \gets \bm{\lambda} (d_t^n, \bm{x})$, $s \gets y-d_r$.
 \EndFor
 \State Compute $c^n = \sum_t \beta^{t-1} c_t^n$. 
\State $n \gets n+1$.
\EndWhile
\State Compute $\tilde{v}^U_T(\bm{e}_i, 0) = \sum_n c^n/nSamples $
\EndFor
\EndProcedure
\end{algorithmic}
\end{algorithm}
\end{minipage}

 Here, $\bm{e}(\mu_0)$ is the unit vector with 1 in the $\mu_0$th position.
 
 \begin{minipage}{\linewidth}
\begin{algorithm}[H]
\begin{algorithmic}[1]
 \caption{Sample Average Cost of Lower Bound Policy $\tilde{v}^L_T(\bm{x}) $}
 \label{algo_lb_tilde}
 \Procedure{$\tilde{v}^L_T$}{$\bm{x}$}
 \For{$i = 1,\dots,N $}  
 \State Set trajectory counter $n \gets 1$. 
 \While {$n<nSamples$ }
\State Initialize $\bm{x} \gets \bm{e}_i$.
 \For {$t=1,\dots,T$} 
\State Calculate $y =  s^*(x)$.
\State Randomly sample the current modulation state $\mu$ with probability $\bm{e}(\mu_0) P_{ij}$. 
 \State Randomly sample the current demand outcome $d^n_t$ with probability $\text{Pr}(d \mid \mu)$.
 \State Compute $c_t^n = h(y-d^n_t)^+ + p(d^n_t-y)^+$.
 \State Update $\mu_0 \gets \mu $, $\bm{x} \gets \bm{\lambda} (d_t^n, \bm{x})$.
 \EndFor
 \State Compute $c^n = \sum_t \beta^{t-1} c_t^n$. 
\State $n \gets n+1$.
\EndWhile
\State Compute $\tilde{v}^L_T(\bm{e}_i) = \sum_n c^n/nSamples $
\EndFor
\EndProcedure
\end{algorithmic}
\end{algorithm}
\end{minipage}

\subsection{Proofs of Results in Section 5}
\begin{proof}[{Proof of \autoref{prop6}}]
The proof of \autoref{prop6} is a direct extension of the results in \cite{scarf60}. 
\end{proof}
\begin{lemma}\label{lemma11}
For all $\bm{x}$ and $n$: 
\begin{description}
\item[(i)] if $s\leq s'$, then $v_n(\bm{x},s) \leq v_n(\bm{x}, s') + K$
\item[(ii)] if $y\leq y'$, then 
$ G_n(\bm{x}, y') - G_n(\bm{x}, y) \geq L(\bm{x}, y') - L(\bm{x}, y) - \beta K $
\item[(iii)] if $s\leq s' \leq \underline{S}(\bm{x})$, then $v_n(\bm{x},s) \geq v_n(\bm{x},s')$
\item[(iv)] if $y\leq y'\leq \underline{S}(\bm{x})$, then
$ G_n(\bm{x},y') - G_n(\bm{x},y) \leq L(\bm{x}, y') - L(\bm{x}, y) \leq 0. $
\end{description}
\end{lemma}
\begin{proof}[{Proof of \autoref{lemma11}.}]
\begin{description}
\item[(i)] This result follows from the $K$-convexity of $v_n(\bm{x},s)$ in $s$, which is a direct implication of the second item of \autoref{prop6}. 
\item[(ii)] This result follows from the definition of $G_n(\bm{x},y)$, the previous result (i), and the fact that $f(y,d)$ is convex and non-decreasing.
\item[(iii)]  $G_n(\bm{x}, s_n(\bm{x})) \leq K+ G_n(\bm{x},S_n(\bm{x})) \leq K+ G_n(\bm{x}, \underline{S}(\bm{x}))$ implies that $s_n(\bm{x}) \leq S_n(\bm{x}) \leq \underline{S}(\bm{x})$ (This is an implication of the definitions of $s_n(\bm{x})$ and $S_n(\bm{x})$, and the fact that $\underline{S}(\bm{x})$ minimizes $L(\bm{x},y)$ while $S_n(\bm{x})$ minimizes the sum of $L(\bm{x},y)$ and a positive term.). It follows from the four cases of $s\leq s' \leq \underline{S}(\bm{x})$ with respect to the value of $s_n(\bm{x})$ that $v_n(\bm{x},s) \geq v_n(\bm{x},s')$.  
\item[(iv)] This result follows from the definition of $G_n(\bm{x},y)$, the non-decreasing nature of $f(y,d)$ in $y$ and  (iii).
\end{description}
\end{proof}
The proof of \autoref{prop7} requires four lemmas. 
\begin{lemma}\label{lemma12}
For all $n$ and $\bm{x}$, $ \underline{S}(\bm{x}) = S_0(\bm{x}) \leq S_n(\bm{x})$. 
\end{lemma}
\begin{lemma}\label{lemma13}
For all $n$ and $\bm{x}$, $s_n(\bm{x})$ can be selected so that $s_n(\bm{x}) \leq \overline{s}(\bm{x})$. 
\end{lemma}
\begin{lemma}\label{lemma14}
For all $n$ and $\bm{x}$, $S_n(\bm{x})$ can be selected so that $S_n(\bm{x}) \leq \overline{S}(\bm{x})$. 
\end{lemma}
\begin{lemma}\label{lemma15}
For all $n$ and $\bm{x}$, $\underline{s}(\bm{x}) \leq s_n(\bm{x})$. 
\end{lemma}
\begin{proof}[{Proof of \autoref{prop7}.}]
The proof of these results follow from the proofs of Lemmas 2 - 5 in \cite{veinott65c}. Proof of \autoref{prop7}(a) follows from Lemmas  \ref{lemma12} -  \ref{lemma15}, and \autoref{prop7}(b) follows from (a) and \autoref{prop6}.  
\end{proof}
\subsection*{Determining {$\Gamma_n(s)$}}
As was true for the $K=0$ case, when $K>0$, there is a finite set of vectors $\Gamma_n(s)$ such that $v_n(\bm{x}, s) = \min \{\bm{x}\bm{\gamma}: \bm{\gamma} \in \Gamma_n(s) \}$ for all $s$. Note that $\Gamma_0(s) = \{\underline{0}\}$ for all $s$, where $\underline{0}$ is the column $N$-vector having zero in all entries. Given $\{ \Gamma_n(s): \forall \ s \}$, we now present an approach for determining $\{\Gamma_{n+1}(s): \forall \ s \}$. Recalling \autoref{sec3.3}, let $\overline{\Gamma} = \{\overline{\bm{\gamma}}_1, \dots, \overline{\bm{\gamma}}_M\}$ be such that $\min_y L(\bm{x},y) = \min\{\bm{x} \bm{\gamma}: \bm{\gamma} \in \overline{\Gamma} \}$. 
Note 
$$ G_n(\bm{x}, y) = L(\bm{x},y) + \beta \sum_{d, z} \sigma(d, z,\bm{x}) v_n \big(\bm{\lambda}(d,z,\bm{x}), f(y,d)\big), $$
for $y \in \{d_1, \dots, d_M\}$. Then, 
$$v_n\big(\bm{\lambda}(d,z,\bm{x}), f(y,d)\big) = \min \{ \bm{\lambda}(d,z,\bm{x}) \bm{\gamma}: \bm{\gamma} \in \Gamma_n(f(y,d)) \}. $$
Let $\Gamma_n'(y)$ be the set of all vectors of the form 
$$ \overline{\bm{\gamma}} + \beta \sum_{d,z} \bm{P}(d,z)\bm{\gamma}(d,z), $$ 
where $\overline{\bm{\gamma}} \in \overline{\Gamma}$ and $\bm{\gamma}(d,z) \in \Gamma_n(f(y,d))$. Then, $G_n(\bm{x},y) = \min\{\bm{x}\bm{\gamma}: \bm{\gamma} \in \Gamma_n'(y)\}$ and 
$$ v_{n+1}(\bm{x},s) = \begin{cases} K+G_n\big(\bm{x},S_n(\bm{x})\big) & \ s\leq s_n(\bm{x}) \\ G_n(\bm{x}, s) & \ \text{otherwise}, \end{cases} $$
where $S_n(\bm{x})$ and $s_n(\bm{x})$ are the smallest integers such that 
\begin{eqnarray*}
G_n \big(\bm{x}, S_n(\bm{x})\big) &\leq& G_n(\bm{x}, y) \ \forall y. 
\\ G_n\big(\bm{x}, s_n(\bm{x})\big) &\leq& K+G_n\big(\bm{x}, S_n(\bm{x})\big).
\end{eqnarray*}
Let $X_n(s', S')$ be the set of all $\bm{x}\in {X}$ such that $s_n(\bm{x}) = s'$ and $S_n(\bm{x}) = S'$. Thus, if $\bm{x}\in X_n(s', S')$, then $s'$ and $S'$ are the smallest integers such that 
\begin{eqnarray*}
G_n \big(\bm{x}, S'(\bm{x})\big) &\leq& G_n(\bm{x}, y) \ \forall y. 
\\ G_n\big(\bm{x}, s'(\bm{x})\big) &\leq& K+G_n\big(\bm{x}, S'(\bm{x})\big).
\end{eqnarray*}
Since $G_n(\bm{x},y)$ is piecewise linear and convex in $\bm{x}$ for each $y$, $X_n(s', S')$ is described by a finite set of linear inequalities. We remark that $\{X_n(s', S'): s' \leq S', \text{ and } X_n(s', S') \neq \emptyset \} $ is a partition of ${X}$. Further, we remark that if $\overline{X}(\underline{s}, \overline{s}, \underline{S}, \overline{S}) \cap X_n(s', S') \neq \emptyset$, then search for $(s', S')$ can be restricted to $\underline{s} \leq s' \leq \overline{s}$ and $\underline{S} \leq S' \leq \overline{S}$. Let $\Gamma_{n+1}(s) = \{ K\underline{1} + \bm{\gamma}: \bm{\gamma} \in \Gamma_n'(S') \}$ for all $s\leq s'$, and let $\Gamma_{n+1}(s) = \Gamma_{n}'(s)$ for all $s>s'$. Thus, $v_{n+1}(\bm{x},s) = \min \{ \bm{x} \bm{\gamma}: \bm{\gamma} \in \Gamma_{n+1}(s) \}$ for all $s$.   

  
  

\end{document}